\documentstyle{amsart}
\font\tenmsb=msbm10
\font\sevenmsb=msbm7
\font\fivemsb=msbm5
\font\largemsb=msbm10 at11pt
\font\largesevenmsb=msbm7 at 7.7pt
\font\largefivemsb=msbm5 at 5.5pt
\font\smallmsb=msbm7 at7.7pt
\newfam\msbfam
\textfont\msbfam=\tenmsb
\scriptfont\msbfam=\sevenmsb
\scriptscriptfont\msbfam=\fivemsb
\newfam\largemsbfam
\textfont\largemsbfam=\largemsb
\scriptfont\largemsbfam=\largesevenmsb
\scriptscriptfont\largemsbfam=\largefivemsb
\newfam\smallmsbfam
\textfont\smallmsbfam=\smallmsb
\def\Bbb#1{\ifdim1em>10.1pt{\fam\largemsbfam\relax#1\kern.9pt}
   \else\ifdim1em<9.9pt{\fam\smallmsbfam\relax#1\kern.5pt}
   \else\fam\msbfam\relax#1\kern.8pt\fi\fi}
\def\qed{\hfill\vbox{\hrule height.09ex
   \hbox{\vrule width.09ex height1.7ex depth.3ex \kern2ex
   \vrule width.09ex height1.7ex depth.3ex}\hrule height.09ex}\bigskip}

\def\Z{{\Bbb Z}}
\def\R{{\Bbb R}}
\def\C{{\Bbb C}}
\def\S{{\Bbb S}}

\newtheorem{theorem}{Theorem}
\newtheorem{proposition}{Proposition}[section]
\newtheorem{lemma}[proposition]{Lemma}
\begin{document}

\large

\title[Sharp Lipschitz estimates for operator $\bar\partial_{\bold M}$]
{Sharp Lipschitz estimates for operator $\bar\partial_{\bold M}$ \\
on a $q$-concave CR manifold}
\author[Peter L. Polyakov]{Peter L. Polyakov}
\address{Department of Mathematics, University of Wyoming, Laramie, WY 82071}
\email{polyakov@@ledaig.uwyo.edu}
\date{October 1996}
\thanks{The author was partially supported by
NSF Grant DMS-9022140 during his stay at MSRI}

\begin{abstract}
We prove that the integral operators $R_r$ and $H_r$ constructed in
\cite{P} and such that
$$f = \bar\partial_{\bold M} R_r(f) +
R_{r+1}(\bar\partial_{\bold M} f) + H_r(f),$$
for a differential form $f \in C_{(0,r)}^{\infty}({\bold M})$
on a regular q-concave CR manifold ${\bold M}$
admit sharp estimates in the Lipschitz scale.
\end{abstract}

\maketitle

\section{Introduction.}\label{Introduction}

\indent
Let ${\bold M}$ be a CR submanifold in a complex
$n$ - dimensional manifold ${\bold G}$ such that for any $z \in {\bold M}$
there exist a neighborhood $V \ni z$ in ${\bold G}$
and smooth real valued functions
$\{ \rho_{k} \} \ \ k = 1, \dots , m  \ ( 1 < m < n-1)$ on $V$ such that
$$\begin{array}{ll}
{\bold M} \cap V = \{ z \in {\bold G} \cap V: \rho_{1}(z) =
\dots = \rho_{m}(z) = 0\},\vspace{0.1in}\\
\partial \rho_1 \wedge \cdots \wedge \partial \rho_m \neq 0 
\ \ \mbox{on} \ {\bold M} \cap V.
\end{array}
\eqno(\arabic{equation})
\newcounter{manifold}
\setcounter{manifold}{\value{equation}}
\addtocounter{equation}{1}$$
\indent
In this paper we continue the study of regularity of the operator
$\bar\partial_{\bold M}$ on a submanifold ${\bold M}$ satisfying special concavity condition. In \cite{P} we proved sharp estimates for solutions
of the $\bar\partial_{\bold M}$ equation with an $L^{\infty}$ right hand
side. Here we push sharp estimates higher on the Lipschitz scale.\\
\indent
Before formulating the main result we will introduce the necessary
notations and definitions.\\
\indent
The CR structure on ${\bold M}$ is induced from ${\bold G}$ and
is defined by the subbundles
$$T^{\prime\prime}({\bold M})
= T^{\prime\prime}({\bold G})|_{{\bold M}} \cap {\bold C}T({\bold M})
\hspace{0.2in} \mbox{and} \hspace{0.2in}
T^{\prime}({\bold M}) = T^{\prime}({\bold G})|_{{\bold M}}
\cap {\bold C}T({\bold M}),$$
\noindent
where ${\bold C}T({\bold M})$ is the complexified tangent bundle
of ${\bold M}$ and the subbundles $T^{\prime\prime}({\bold G})$ and
$T^{\prime}({\bold G})= \overline{T^{\prime\prime}}({\bold G})$ of
the complexified tangent bundle ${\bold C}T({\bold G})$ define the complex
structure on ${\bold G}$.\\
\indent
We will denote by $T^c({\bold M})$ the subbundle $T({\bold M}) \cap
\left[ T^{\prime}({\bold M}) \oplus T^{\prime\prime}({\bold M})\right].$
If we fix a hermitian scalar product on ${\bold G}$ then we can choose
a subbundle $N \in T({\bold M})$ of real dimension $m$ such that
$T^c({\bold M}) \perp N$ and for a complex subbundle 
${\bold N} = {\bold C}N$ of ${\bold C}T({{\bold M}})$ we have
$${\bold C}T({{\bold M}}) =
T^{\prime}({\bold M}) \oplus T^{\prime\prime}({\bold M})
\oplus {\bold N},
\hspace{0.1in}T^{\prime}({\bold M})
\perp {\bold N} \hspace{0.1in}\mbox{and}
\hspace{0.1in}T^{\prime\prime}({\bold M}) \perp {\bold N}.$$
\indent
We define the Levi form of ${\bold M}$ as the hermitian form on
$T^{\prime}({\bold M})$ with values in ${\bold N}$
$${\cal L}_z(L(z)) = \sqrt{-1} \cdot \pi
\left( \left[\overline L, L\right]\right)(z)
\hspace{0.2in}\left( L(z) \in T^{\prime}_z({\bold M})\right),$$
\noindent
where $\left[ \overline L, L \right] =
\overline L L - L \overline L$ and $\pi$ is the projection
of ${\bold C}T({\bold M})$ along
$T^{\prime}({\bold M}) \oplus T^{\prime\prime}({\bold M})$
onto $N$.\\
\indent
If the functions $\{ \rho_{k} \}$ are chosen so that the vectors
$\{ \mbox{grad} \rho_{k} \}$ are orthonormal then 
the Levi form of ${\bold M}$ may be defined as
$$L_z({\bold M}) = - \sum_{k=1}^m\left(L_{z}\rho_k(\zeta)\right) \cdot 
\mbox{grad $\rho_k(z)$},$$
\noindent
where $L_{z}\rho(\zeta)$ is the Levi form
of the real valued function $\rho \in C^4({\bold D})$ at the point $z$:
$$L_{z}\rho(\zeta) = \sum_{i,j}\frac{{\partial}^2\rho}
{\partial \zeta_{i} \partial\bar \zeta_{j}}(z) \ \zeta_{i} 
\cdot \bar{\zeta_{j}}.$$
\indent
Analogously by the Hessian form of ${\bold M}$ at the point
$z \in {\bold M}$ we call the
hermitian form on the complex tangent space $T^c_{z}({\bold M})$ 
of ${\bold M}$ at $z$ with values in $N_{z}$, defined by the formula:
$${\cal H}_z({\bold M}) =
- \sum_{k=1}^m \left( {\cal H}_{z}\rho_k(\zeta) \right) \cdot
\mbox{grad $\rho_k(z)$},$$
\noindent
where ${\cal H}_{z}\rho(\zeta)$ is the Hessian form
of the real valued function $\rho \in C^4({\bold D})$ at the point $z$:
$${\cal H}_{z}\rho(\zeta) = \sum_{i,j}\frac{{\partial}^2\rho}
{\partial \zeta_{i} \partial\bar \zeta_{j}}(z) \ \zeta_{i} 
\cdot \bar{\zeta_{j}} + \mbox{Re}\left\{ \sum_{i,j}\frac{{\partial}^2\rho}
{\partial \zeta_{i} \partial \zeta_{j}}(z) \ \zeta_{i} 
\cdot {\zeta_{j}} \right\}.$$
\indent
For a pair of vectors $\mu = (\mu_1, \dots, \mu_n )$ and
$\nu = (\nu_1, \dots, \nu_n )$ we will denote
$\langle \mu, \nu \rangle = \sum_{i=1}^n \mu_i \cdot \nu_i$.\\
\indent
For a unit vector
$\theta = (\theta_1, \dots \theta_m) \in \mbox{Re}{\bold N}_z$ we
define the Levi form of 
${\bold M}$ at the point $z \in {\bold M}$ in the direction 
$\theta$ as the scalar hermitian form on
$T^c_{z}({\bold M})$
$$\langle \theta,\ L_z({\bold M}) \rangle = - L_{z}\rho_\theta(\zeta),$$
and the Hessian form of 
${\bold M}$ at the point $z \in {\bold M}$ in the direction 
$\theta$ as the hermitian form
$$\langle \theta,\ {\cal H}_z({\bold M}) \rangle =
- {\cal H}_{z}\rho_\theta(\zeta),$$
\noindent
where $\rho_\theta(\zeta) = \sum_{k=1}^m \theta_k \rho_k(\zeta)$.\\
\indent
Following \cite{H2} we introduce the notions of q-pseudoconcave and q-concave
CR manifolds.\\
\indent
Namely, we call ${\bold M}$ q-pseudoconcave
(weakly q-pseudoconcave) at $z \in {\bold M}$ in the 
direction $\theta$ if the Levi form of ${\bold M}$ at $z$ in this direction 
$\langle \theta,\ L_z({\bold M}) \rangle$ has at least $q$ negative 
($q$ nonpositive) eigenvalues on $T^c_{z}({\bold M})$ and
we call ${\bold M}$ q-concave (weakly q-concave) at $z \in {\bold M}$ in the 
direction $\theta$ if the Hessian form of ${\bold M}$ at $z$ in this direction
$\langle \theta,\ {\cal H}_z({\bold M}) \rangle$ has at least $q$ negative 
($q$ nonpositive) eigenvalues on $T^c_{z}({\bold M})$.\\
\indent
We call ${\bold M}$ q-pseudoconcave (weakly q-pseudoconcave)
at $z \in {\bold M}$ if it is q-pseudoconcave (weakly q-pseudoconcave)
in all directions. Analogously, we call ${\bold M}$ q-concave
(weakly q-concave) at $z \in {\bold M}$ if it is q-concave
(weakly q-concave) in all directions.\\
\indent
We call a q-concave CR manifold ${\bold M}$ by a
regular q-concave CR manifold (cf. \cite{P})
if for any $z \in {\bold M}$ there exist an open neighborhood $U \ni z$
in  ${\bold M}$ and a family $E_{q}(\theta,z)$ of
$q$-dimensional complex linear subspaces in $T^c_{z}({\bold M})$
smoothly depending on $(\theta, z) \in \S^{m-1} \times U$ and such that
the Hessian form $\langle \theta,\ {\cal H}_z({\bold M}) \rangle$
is strictly negative on $E_{q}(\theta,z).$\\
\indent
Following \cite{S} we define spaces
$\Gamma^{p, \alpha}({\bold M})$
for nonnegative $p \in \Z$ and $0 < \alpha <2$. Namely, we say that function
$h \in \Gamma^{p,\alpha}({\bold M})$ if for any set of tangent vector fields
$D_1, \dots , D_p$ on ${\bold M}$ such that
$\| D_i \|_{ C^{p+2}\left({\bold M}\right) } \leq 1$
$$\| h \|_{ \Gamma^{p, \alpha}\left({\bold M}\right) } =
\| D_1 \circ \cdots \circ D_p h \|_{\Lambda^{ \frac{\alpha}{2} }
\left({\bold M}\right) } + \sup \left\{
\| D_1 \circ \cdots \circ D_p h(x(\cdot))\|_{\Lambda^{\alpha}([0,1])} \right\}
< \infty,$$
\noindent
where the {\it sup} is taken over all curves $x : [0,1] \rightarrow
{\bold M}$ such that
$$\begin{array}{ll}
(i)\hspace{0.1in} |x^{\prime}(t)|, |x^{\prime\prime}(t)| \leq 1,
\vspace{0.1in}\\
(ii)\hspace{0.1in}x^{\prime}(t) \in T^{c}({\bold M}).
\end{array}$$
\indent
For a differential form $g = \sum_{I,J}g_{I,J}(z)dz^I \wedge d\bar{z}^J$
with $|I| = k$ and $|J| = r$ we say that
$g \in \Gamma^{p,\alpha}_{(k,r)}({\bold M})$ if $g_{I,J} \in
\Gamma^{p,\alpha}({\bold M}).$\\
\indent
The following theorem represents the main result of the paper.\\

\begin{theorem}\label{HomotopyTheorem}
Let $0< \alpha <1$ and let a $C^{\infty}$ submanifold
${\bold M} \subset {\bold G}$
of the form (\arabic{manifold}) be regular q-concave with $q \geq 2$.
Then for any $r = 1,\dots, q-1$ there exist linear operators
$${\bold R}_r: \Gamma^{p, \alpha}_{(0,r)}({\bold M}) \rightarrow
\Gamma^{p, \alpha +1}_{(0,r-1)}({\bold M})\hspace{0.1in}
\mbox{and} \hspace{0.1in}
{\bold H}_r: \Gamma^{p, \alpha}_{(0,r)}({\bold M}) \rightarrow
\Gamma^{p, \alpha}_{(0,r)}({\bold M})$$\\
such that ${\bold R}_r$ is bounded and ${\bold H}_r$ is compact and such that
for any differential form $f \in
C_{(0,r)}^{\infty}({\bold M})$ the equality:
$$f = \bar\partial_{\bold M} {\bold R}_r(f) +
{\bold R}_{r+1}(\bar\partial_{\bold M} f) + {\bold H}_r(f)
\eqno(\arabic{equation})
\newcounter{HomotopyFormula}
\setcounter{HomotopyFormula}{\value{equation}}
\addtocounter{equation}{1}$$
holds.
\end{theorem}

\indent
The study of the $\bar\partial_{\bold M}$ complex on a real submanifold
${\bold M} \in \C^n$ was initiated by J. J. Kohn
and H. Rossi in \cite{K}, \cite{KR}. For a closed strongly
pseudoconvex hypersurface ${\bold M} \in \C^n$ J. J. Kohn \cite{K},
A. Andreotti, C. D. Hill \cite{AnH}, I. Naruki \cite{Na}
proved the
solvability of the nonhomogeneous $\bar\partial_{\bold M}$ equation:
$$\bar\partial_{\bold M} f^{0,r-1} = g^{0,r}$$\\
for $r < n-1$ and any $g^{0,r} \in C_{(0,r)}^{\infty}({\bold M})$
satisfying the solvability condition
$$\bar\partial_{\bold M} g^{0,r} = 0.$$\\
\indent
A version of the main theorem for q-pseudoconvex hypersurfaces $(m=1)$
and different spaces $\Gamma^{\alpha}({\bold M})$ with $\alpha \geq 0$
was proved by G. B. Folland and E. M. Stein in \cite{FS}
(cf. also \cite{H1}).\\
\indent
Generalizing notion of q-concavity to
the manifolds of higher codimension I. Naruki in
\cite{Na} using Kohn-Hormander's method constructed bounded operators
$R_r: L_{(0,r)}^{2}({\bold M})
\rightarrow L_{(0,r-1)}^{2}({\bold M^\prime})$
for the $(0,r)$ forms with $r > n-m-q$.\\
\indent
Then in \cite{H2} and \cite{AH} with the use of explicit integral formulas
bounded operators $R_r: L_{(0,r)}^{\infty}({\bold M})
\rightarrow \Gamma_{(0,r-1)}^{0,1-\epsilon}({\bold M})$ were
constructed on a q-pseudoconcave CR manifold of higher codimension
for the forms of type $(0,r)$ with $r < q$ or $r > n-m-q$.\\
\indent
Existence of solution $f \in C_{(0,r-1)}^{\frac{1}{2}}({\bold M})$ of
equation $\bar\partial_{\bold M} f = g$ for
$g \in \L_{(0,r)}^{\infty}({\bold M})$ for q-pseudoconcave CR manifold
was obtained by M. Y. Barkatou in \cite{Ba} and existence of solution
$f \in \Gamma_{(0,r-1)}^{\beta+1}({\bold M})$ for
$g \in \Gamma_{(0,r)}^{\beta}({\bold M})$, $0 < \beta < 1$ and ${\bold M}$ -
quadratic q-pseudoconcave CR manifold was obtained in the
paper \cite{BGG} by R. Beals, B. Gaveau and P.C. Greiner.\\
\indent
The problem considered here is closely related to the problem of obtaining
sharp Lipschitz estimates for solutions of $\bar\partial$ equation on
complex manifolds with boundary. Paper by Y.-T. Siu \cite{Siu}
devoted to this problem has motivated some of the constructions
in our paper.\\
\indent
Author thanks R. Beals and G. Henkin for helpful discussions.\\

\section{Construction of ${\bold R}_r$ and ${\bold H}_r$.}
\label{Construction}

\indent
The local versions of operators ${\bold R}_r$ and ${\bold H}_r$ from the
theorem~\ref{HomotopyTheorem} were constructed in \cite{P}. Before describing
these local operators we will introduce necessary definitions and notations.\\
\indent
For a vector-valued function $\eta =(\eta_1,\dots,\eta_n)$ we will use
the notation:
$$\omega^{\prime}(\eta) = \sum_{k=1}^{n} (-1)^{k-1} \eta_k
\wedge_{j \neq k} d\eta_j, \hspace{0.2in} \omega(\eta) =
\wedge_{j = 1}^n d\eta_j.$$
\indent
If $\eta = \eta(\zeta, z,t)$ is a smooth function of $\zeta \in \C^n,
z \in \C^n$ and a real parameter $t \in \R^p$ satisfying
the condition 
$$\sum_{k=1}^{n} \eta_k(\zeta, z,t) \cdot (\zeta_k - z_k) = 1
\eqno(\arabic{equation})
\newcounter{leray}
\setcounter{leray}{\value{equation}}
\addtocounter{equation}{1}$$
then
$$d \omega^{\prime}(\eta) \wedge \omega(\zeta) \wedge \omega(z) = 0$$\\
or, separating differentials,
$$d_{t} \omega^{\prime}(\eta) + \bar\partial_{\zeta} 
\omega^{\prime}(\eta) + \bar\partial_z \omega^{\prime}(\eta) = 0.
\eqno(\arabic{equation})
\newcounter{dofomega}
\setcounter{dofomega}{\value{equation}}
\addtocounter{equation}{1}$$
\indent
Also, if $\eta(\zeta,z,t)$ satisfies (\arabic{leray}) then the
differential form $\omega^{\prime}(\eta) \wedge \omega(\zeta) 
\wedge \omega(z)$ can be represented as:
$$\sum_{r=0}^{n-1} \omega^{\prime}_r(\eta) \wedge \omega(\zeta)
\wedge \omega(z),
\eqno(\arabic{equation})
\newcounter{sumofomega}
\setcounter{sumofomega}{\value{equation}}
\addtocounter{equation}{1}$$
where $\omega^{\prime}_r(\eta)$ is a differential form of the order
$r$ in $d \bar{z}$ and respectively of the order $n-r-1$ in 
$d \bar\zeta$ and $dt$. From (\arabic{dofomega}) and 
(\arabic{sumofomega}) follow equalities:
$$ d_t \omega^{\prime}_r(\eta) + \bar\partial_{\zeta} 
\omega^{\prime}_r(\eta) + \bar\partial_z \omega^{\prime}_{r-1}(\eta) = 0
\hspace{0.3in} (r=1,\dots,n),
\eqno(\arabic{equation})
\newcounter{domega}
\setcounter{domega}{\value{equation}}
\addtocounter{equation}{1}$$
and
$$\omega^{\prime}_r(\eta) = \frac{1}{(n-r-1)!r!}
\mbox{Det} \left[\eta, \hspace{0.05in}
\overbrace{\bar\partial_z \eta}^{r}, \hspace{0.05in} 
\overbrace{\bar\partial_{\zeta,t} \eta}^{n-r-1}
\right],
\eqno(\arabic{equation})
\newcounter{determinantomega}
\setcounter{determinantomega}{\value{equation}}
\addtocounter{equation}{1}$$
where the determinant is calculated by the usual rules but with external
products of elements and the position of the element in the external product
is defined by the number of its column.\\
\indent
Let ${\widetilde U}$ be an open neighborhood in ${\bold G}$ and
$U = {\widetilde U} \cap {\bold M}$. We call a vector function
$$P(\zeta,z) = (P_1(\zeta,z),\dots, P_n(\zeta,z)) \hspace{0.2in}
\mbox{for} \hspace{0.2in} (\zeta,z) \in
\left( {\widetilde U} \setminus U \right) \times {\widetilde U}$$
by strong ${\bold M}$-barrier for ${\widetilde U}$ if there exist
$C,c >0$ such that the inequality:
$$\left| \Phi(\zeta,z) \right| > C \cdot \left( \rho(\zeta)
+ | \zeta - z|^2 \right)^2,
\eqno(\arabic{equation})
\newcounter{barrier}
\setcounter{barrier}{\value{equation}}
\addtocounter{equation}{1}$$
holds for $(\zeta,z) \in \left({\widetilde U} \setminus U \right)
\times U$, where
$$\Phi(\zeta,z) = \langle P(\zeta,z), \zeta - z \rangle = 
\sum_{i=1}^n P_i(\zeta,z) \cdot (\zeta_i - z_i).$$\\
\indent
In \cite{P} special strong ${\bold M}$-barriers for a q-concave
CR submanifold of the form (1) were constructed. Below we describe
these barriers.\\
\indent
According to (\arabic{manifold}) we may assume that
$U = {\widetilde U} \cap {\bold M}$ is a
set of common zeros of smooth functions
$\{ \rho_k, \hspace{0.1in}k=1, \dots, m \}$.
Then, using the q-concavity of ${\bold M}$ and applying
Kohn's lemma to the set of functions $\{ \rho_k \}$ we can
construct a new set of functions $\tilde \rho_1, \dots , \tilde
\rho_m$ of the form:
$$\tilde \rho_k(z) = \rho_k(z) + A \cdot
\left( \sum_{i=1}^m \rho_i^2(z) \right),$$\\ 
with large enough constant $A > 0$ and such that
for any $z \in {\bold M}$ there exist an open neighborhood
$U \ni z$ and a family $E_{q+m}(\theta,z)$ of
$q+m$ dimensional complex linear subspaces in $\C^n$
smoothly depending on $(\theta, z) \in \S^{m-1} \times U$ and such that
$- {\cal H}_{z}\tilde\rho_{\theta}$ is strictly negative
on $E_{q+m}(\theta,z)$ with all negative eigenvalues not exceeding some
$c < 0$.\\
\indent
To simplify notations we will assume that the functions
$\rho_1, \dots , \rho_m$ already satisfy this condition.\\
\indent
Let $E_{n-q-m}^{\bot}(\theta,z)$ be the family of $n-q-m$ dimensional
subspaces in $T({\bold G})$ orthogonal to $E_{q+m}(\theta,z)$ and let
$$a_j(\theta,z) = \left( a_{j1}(\theta,z), \dots, a_{jn}(\theta,z)
\right) \hspace{0.05in} \mbox{for} \hspace{0.05in}j = 1, \dots, n-q-m$$
be a set of $C^2$ smooth vector functions
representing an orthonormal basis in $E_{n-q-m}^{\bot}(\theta,z).$\\
\indent
Defining for $(\theta, z, w) \in \S^{m-1} \times U \times \C^n$
$$A_j(\theta,z,w) = \sum_{i=1}^n a_{ji}(\theta,z) \cdot w_i,
\hspace{0.1in}(j=1, \dots, n-q-m)$$
we construct the form
$${\cal A}(\theta, z, w) = \sum_{j=1}^{n-q-m} A_j(\theta,z,w) \cdot
{\bar A_j(\theta,z,w)}$$
such that the hermitian form
$${\cal H}_{z}\rho_{\theta}(w) + {\cal A}(\theta, z, w)$$
is strictly positive definite in $w$ for $(\theta,z) \in \S^{m-1} \times U$.\\
\indent
Then we define for $\zeta, z \in U$:
$$\begin{array}{llllll}
\hspace{0.4in}Q^{(k)}_i(z) = - \frac{\partial \rho_k}
{\partial \zeta_i}(z),\vspace{0.2in}\\
\hspace{0.4in}F^{(k)} (\zeta, z) = \langle Q^{(k)}(z), \zeta - z \rangle,
\vspace{0.2in}\\
\hspace{0.4in}P_i(\zeta, z) = \sum_{k=1}^m
Q^{(k)}_i(z) \cdot {\overline F}^{(k)} (\zeta, z)\vspace{0.2in}\\
+ \sum_{j=1}^{n-q-m}
a_{ji}(\theta(\zeta), z) \cdot {\bar A_j(\theta(\zeta),z,\zeta-z)}
\cdot {\cal A}(\theta(\zeta),z,\zeta-z),\vspace{0.2in}\\
\hspace{0.4in}\Phi (\zeta, z) = \langle P(\zeta, z), 
\zeta - z \rangle\vspace{0.2in}\\
= \sum_{k=1}^m F^{(k)} (\zeta, z)
\cdot {\overline F}^{(k)} (\zeta, z) +
\frac{1}{4} {\cal A}^2(\theta(\zeta),z,\zeta-z)
\end{array}
\eqno(\arabic{equation})
\newcounter{Barrier}
\setcounter{Barrier}{\value{equation}}
\addtocounter{equation}{1}$$
with 
$$\theta_k(\zeta) = - \frac{\rho_k(\zeta)}{\rho(\zeta)} 
\hspace{0.2in} \mbox{for $k=1,\dots, m$}.$$
\indent
From the definitions of $Q_i^{(s)}(z)$, $F^{(s)}(\zeta,z)$ and
${\cal A}(\zeta,z)$ follow the equalities that will
be used in the further estimates
$$\sum_{i=1}^n Q_i^{(s)}(z)d\zeta_i = d_{\zeta}F^{(s)}(\zeta, z),
\eqno(\arabic{equation})
\newcounter{dzetaF}
\setcounter{dzetaF}{\value{equation}}
\addtocounter{equation}{1}$$
$$\bar\partial_{\zeta}({\cal A}{\bar A_j})=
\mu_{\tau}^{(j)}(\zeta,z) + \mu_{\nu}^{(j)}(\zeta,z)
\eqno(\arabic{equation})
\newcounter{mu}
\setcounter{mu}{\value{equation}}
\addtocounter{equation}{1}$$
where
$$\mu_{\tau}^{(j)}(\zeta,z)= {\cal A}(\zeta,z) \cdot
\sum_{i=1}^n {\bar a_{ji}}(\theta(\zeta),z) d{\bar \zeta_i}+$$
$$+ {\bar A_j}(\zeta,z) \cdot \sum_{i=1}^n \sum_{k=1}^{n-q-m}
A_k(\zeta,z) {\bar a_{ki}}(\theta(\zeta),z) d{\bar \zeta_i},$$
$$\mu_{\nu}^{(j)}(\zeta,z)= {\cal A}(\zeta,z) \cdot \sum_{i=1}^n
(\bar\zeta_i - \bar z_i) \bar\partial_{\zeta} {\bar a_{ji}}
(\theta(\zeta),z)+$$
$$+{\bar A_j}(\zeta,z) \cdot \sum_{i=1}^n (\bar\zeta_i - \bar z_i)
\sum_{k=1}^{n-q-m} A_k(\zeta,z) \bar \partial_{\zeta} \bar a_{ki}
(\theta(\zeta),z)+$$
$$+{\bar A_j}(\zeta,z) \cdot \sum_{i=1}^n (\zeta_i - z_i)
\sum_{k=1}^{n-q-m} \bar A_k(\zeta,z) \bar \partial_{\zeta} a_{ki}
(\theta(\zeta),z),$$
\noindent
and
$$\left\{\begin{array}{ll}
\bar \partial_{\zeta} {\overline F}^{(s)}(\zeta,z)
= -\bar\partial_{\zeta} \rho_s(\zeta) + \chi_s(\zeta,z),\vspace{0.1in}\\
\bar\partial_{z} {\overline F}^{(s)}(\zeta,z)
= \bar\partial_{z} \rho_s(z) + \kappa_s(\zeta,z),
\end{array} \right.
\eqno(\arabic{equation})
\newcounter{dbarFbar}
\setcounter{dbarFbar}{\value{equation}}
\addtocounter{equation}{1}$$
where
$$\chi_s(\zeta,z) = \sum_{i=1}^n ({\bar \zeta}_i -
{\bar z}_i) \cdot \chi_{is}(\zeta,z)$$\\
are differential forms in $\zeta$ with smooth coefficients and
$$\kappa_s(\zeta,z) = \sum_{i=1}^n ({\bar \zeta}_i -
{\bar z}_i) \cdot \kappa_{is}(\zeta,z)$$\\
are differential forms in $z$ with smooth coefficients.\\
\indent
In our construction of the global integral formula on ${\bold M}$ we will use
local formula from \cite{P}. To describe this formula we introduce the
following notations.\\
\indent
We define the tubular neighborhood ${\bold G_{\epsilon}}$ of ${\bold M}$ in 
${\bold G}$ as follows:
$${\bold G}_{\epsilon}= \{z \in  {\bold G}: \rho(z) < \epsilon \},$$
\noindent
where $\rho(z) = {\left(\sum_{k=1}^m \rho_k^2(z)\right)}^{\frac{1}{2}}$.
The boundary of ${\bold G}_{\epsilon}$ -- ${\bold M}_{\epsilon}$ is defined
by the condition
$${\bold M}_\epsilon = \{ z \in {\bold G}: \rho(z) = \epsilon \}.$$
\indent
For a sufficiently small neighborhood
${\widetilde U} \in {\bold G}$ we may assume that functions
$$\rho_k(\zeta),\hspace{0.05in}\mbox{Im}F^{(k)}(\zeta,z) \hspace{0.05in}
\{ k=1,\dots, m \}$$
have a nonzero jacobian in $\mbox{Re}\zeta_{i_1}, \dots,
\mbox{Re}\zeta_{i_m},\hspace{0.05in}\mbox{Im}\zeta_{i_1}, \dots,
\mbox{Im}\zeta_{i_m}$ for $z, \zeta \in {\widetilde U}$.
Therefore, for any fixed $z \in {\widetilde U}$ these functions 
may be chosen as local $C^{\infty}$ coordinates in $\zeta$.
We may also complement the functions above by holomorphic
functions $w_j(\zeta)=u_j(\zeta) + iv_j(\zeta)$ with $j=1,\dots, n-m$
so that the functions
$$\rho_k(\zeta),\hspace{0.05in}\mbox{Im}F^{(k)}(\zeta,z)
\hspace{0.05in} \{ k=1,\dots, m \},$$
$$u_j(\zeta), v_j(\zeta)\hspace{0.05in} \{ j=1,\dots, n-m \},$$
represent a complete system of local coordinates in
$\zeta \in {\widetilde U}$ for any fixed $z \in {\widetilde U}$.\\
\indent
The following complex valued vector fields on $U$ for any fixed $z \in U$
$$Y_{i,\zeta}(z) = \frac{\partial}{\partial\mbox{Im}F^{(i)}(\zeta,z)}
\hspace{0.05in}\mbox{for}\hspace{0.05in}i=1,\dots, m,$$
$$\frac{\partial}{\partial w_i},\hspace{0.05in}
\frac{\partial}{\partial {\bar w}_i}\hspace{0.05in}\mbox{for}\hspace{0.05in}
i=1,\dots, n-m,$$
represent a basis in ${\bold C}T({\bold M})$ at any $z \in U$.\\
\indent
We need also to consider local extensions
of functions and forms from $U = {\widetilde U} \cap {\bold M}$ to
${\widetilde U}$. Defining appropriate functional
spaces on these neighborhoods we say that a function $h$ on
${\widetilde U}(\epsilon) = {\bold G}_{\epsilon} \cap {\widetilde U}$
is in $\Gamma^{p, \alpha} \left( \{ \rho \},
{\widetilde U}(\epsilon) \right)$ if
$$\left\| h \right\|_{ \Gamma^{p, \alpha} \left( \{ \rho \},
{\widetilde U}(\epsilon) \right) }
= \left\| h \right\|_{\Lambda^{p+\alpha/2}({\widetilde U}(\epsilon))}
+ \sup_{ \left( \sum_{k=1}^m \delta_k^2 \right)^{1/2} < \epsilon }
\left\| h \right\|_{\Gamma^{p, \alpha}
(U(\delta_1,\dots,\delta_m))} < \infty$$
\noindent
where
$$U(\delta_1,\dots,\delta_m) = \left\{ z \in {\widetilde U}
: \rho_1(z)=\delta_1,\dots,\rho_m(z)= \delta_m \right\}.$$
\indent
For a differential form $f = \sum_{I,J}f_{I,J}(z)dz^I \wedge d\bar{z}^J$
with $|I| = l$ and $|J| = r$ we say that
$f \in \Gamma^{p, \alpha}_{(l,r)} \left( \{ \rho \},
{\widetilde U}(\epsilon) \right)$ if $f_{I,J} \in
\Gamma^{p, \alpha} \left( \{ \rho \}, {\widetilde U}(\epsilon) \right).$\\
\indent
We introduce a local extension operator for
${\widetilde U}(\epsilon)$ and $U = {\widetilde U}(\epsilon) \cap {\bold M}$
$$E_{U} : \Gamma^{p, \alpha}_{(0,r)} (U)
\rightarrow \Gamma^{p, \alpha}_{(0,r)}
\left( \{ \rho \}, {\widetilde U}(\epsilon) \right)$$
which we define by extending all the coefficients of the differential form
identically with respect to coordinates, complementary to
$\rho_1, \dots, \rho_m$ in ${\widetilde U}$. From the construction it follows
that $E_{U}$ satisfies the following estimate
$$\| E_{U}(g) \|_{\Gamma^{p, \alpha}
\left( \{ \rho \}, {\widetilde U}(\epsilon) \right) } \leq
C \cdot \| g \|_{ \Gamma^{p, \alpha}(U) }.$$
\indent
The following proposition provides local integral formula for
$\bar\partial_{\bold M}$.

\begin{proposition}\label{LocalHomotopy}\hspace{0.1in}\cite{P}\hspace{0.1in}
Let ${\bold M} \subset {\bold G}$ be a $C^{\infty}$ regular q-concave
CR submanifold of the form (\arabic{manifold}) and ${\widetilde U}$ an open
neighborhood in ${\bold G}$ with analytic coordinates $z_1, \dots, z_n$.\\
\indent
Then for $r = 1,...,q-1$ and any differential form
$g \in C_{(0,r)}^{\infty}({\bold M})$ with compact support in
$U$ the following equality
$$g = \bar\partial_{\bold M} R_r(g) +
R_{r+1}(\bar\partial_{\bold M} g) + H_r(g),
\eqno(\arabic{equation})
\newcounter{LocalFormula}
\setcounter{LocalFormula}{\value{equation}}
\addtocounter{equation}{1}$$
holds, where
$$R_r(g)(z)$$
$$= (-1)^r \frac{(n-1)!}{(2\pi i)^n} \cdot \mbox{pr}_{\bold M} \circ
\lim_{\epsilon \rightarrow 0}
\int_{{\bold M}_{\epsilon}\times [0,1]} {\widetilde g}(\zeta)
\wedge\omega^{\prime}_{r-1}\left((1-t)\frac{\bar\zeta - \bar z}
{{\mid \zeta - z \mid}^2} + t\frac{P(\zeta,z)}
{\Phi(\zeta,z)}\right) \wedge\omega(\zeta),$$
$$H_r(g)(z) = (-1)^{r} \frac{(n-1)!}{(2\pi i)^n}
\cdot \mbox{pr}_{\bold M} \circ \lim_{\epsilon \rightarrow 0}
\int_{{\bold M}_{\epsilon}} {\widetilde g}(\zeta)
\wedge \omega^{\prime}_{r} \left( \frac{P(\zeta,z)}
{\Phi(\zeta,z)}\right) \wedge\omega(\zeta),$$\\
\noindent
${\widetilde g} = E_{U}(g)$ is the extension of $g$,
$\Phi(\zeta,z)$ is a local barrier for ${\widetilde U}$ constructed
in (\arabic{Barrier}) and $\mbox{pr}_{\bold M}$ denotes the operator
of projection to the space of tangential differential forms on ${\bold M}$.
\end{proposition}

\indent
{\bf Remark.} As follows from the definitions above spaces
$\Gamma^{p, \alpha} \left( \{ \rho \},{\widetilde U}(\epsilon) \right)$
as well as the extension operator $E_{U}$ depend on the choice of
functions $\rho_1, \dots, \rho_m$. But we notice (cf. \cite{P}) that
the operators $R_r$ and $H_r$ are independent of the choice of
functions $\rho_1, \dots, \rho_m$ and of extension operator $E_{U}$.\\

\indent
To construct now global formula on ${\bold M}$ we consider two finite
coverings
$\{ {\widetilde U}_{\iota} \subset {\widetilde U}^{\prime}_{\iota}\}$ of
${\bold G}$ and two partitions of unity $\{\vartheta_{\iota} \}$ and
$\{\vartheta^{\prime}_{\iota} \}$ subordinate to these coverings and
such that $\vartheta^{\prime}_{\iota}(z) = 1$ for
$z \in \mbox{supp}(\vartheta_{\iota})$.\\
\indent
Applying proposition~\ref{LocalHomotopy} to the form $\vartheta_{\iota}g$ in
$U^{\prime}_{\iota}$ we obtain
$$\vartheta_{\iota}(z)g(z) =
\bar\partial_{\bold M} R_r^{\iota}(\vartheta_{\iota}g)(z)
+ R_{r+1}^{\iota}(\bar\partial_{\bold M} \vartheta_{\iota}g)(z)
+ H_r^{\iota}(\vartheta_{\iota}g)(z).$$
\indent
Multiplying the equality above by $\vartheta^{\prime}_{\iota}(z)$ and
using equalities
$$\vartheta^{\prime}_{\iota}(z) \cdot \bar\partial_{\bold M}
R_r^{\iota}(\vartheta_{\iota}g)(z) = \bar\partial_{\bold M} \left[
\vartheta^{\prime}_{\iota}(z) \cdot
R_r^{\iota}(\vartheta_{\iota}g)(z) \right] - \bar\partial_{\bold M}
\vartheta^{\prime}_{\iota}(z) \wedge R_r^{\iota}(\vartheta_{\iota}g)(z)$$
\noindent
and
$$R_{r+1}^{\iota}(\bar\partial_{\bold M} \vartheta_{\iota}g)(z) =
R_{r+1}^{\iota}(\bar\partial_{\bold M} \vartheta_{\iota} \wedge g)(z)
+ R_{r+1}^{\iota}(\vartheta_{\iota} \bar\partial_{\bold M} g)(z)$$
we obtain
$$\vartheta_{\iota}(z)g(z) = \bar\partial_{\bold M}
{\bold R}_r^{\iota}(g)(z)+{\bold R}_{r+1}^{\iota}(\bar\partial_{\bold M}g)(z)
+{\bold H}_r^{\iota}(g)(z)
\eqno(\arabic{equation})
\newcounter{varthetag}
\setcounter{varthetag}{\value{equation}}
\addtocounter{equation}{1}$$
\noindent
with
$${\bold R}_r^{\iota}(g)(z) = \vartheta^{\prime}_{\iota}(z) \cdot
R_r^{\iota}(\vartheta_{\iota}g)(z)$$
\noindent
and
$${\bold H}_r^{\iota}(g)(z) = - \bar\partial_{\bold M}
\vartheta^{\prime}_{\iota}(z) \wedge R_r^{\iota}(\vartheta_{\iota}g)(z)
+ \vartheta^{\prime}_{\iota}(z) \cdot
R_{r+1}^{\iota}(\bar\partial_{\bold M} \vartheta_{\iota} \wedge g)(z)
+ \vartheta^{\prime}_{\iota}(z) \cdot H_r^{\iota}(\vartheta_{\iota}g)(z).$$
\indent
Adding equalities (\arabic{varthetag}) for all $\iota$ we obtain

\begin{proposition}\label{GlobalHomotopy}
Let ${\bold M} \subset {\bold G}$ be a $C^{\infty}$ regular q-concave
CR submanifold of the form (\arabic{manifold}).\\
\indent
Then for $r = 1,...,q-1$ and any differential form
$g \in C_{(0,r)}^{\infty}({\bold M})$ the following equality
$$g = \bar\partial_{\bold M} {\bold R}_r(g) +
{\bold R}_{r+1}(\bar\partial_{\bold M} g) + {\bold H}_r(g),
\eqno(\arabic{equation})
\newcounter{GlobalFormula}
\setcounter{GlobalFormula}{\value{equation}}
\addtocounter{equation}{1}$$
holds, where
$${\bold R}_r(g)(z) = \sum_{\iota} \vartheta^{\prime}_{\iota}(z) \cdot
R_r^{\iota}(\vartheta_{\iota}g)(z)$$
\noindent
and
$${\bold H}_r(g)(z) = \sum_{\iota} \left[ - \bar\partial_{\bold M}
\vartheta^{\prime}_{\iota}(z) \wedge R_r^{\iota}(\vartheta_{\iota}g)(z)
+ R_{r+1}^{\iota}(\bar\partial_{\bold M} \vartheta_{\iota} \wedge g)(z)
+ \vartheta^{\prime}_{\iota}(z) \cdot H_r^{\iota}(\vartheta_{\iota}g)(z)
\right].$$
\end{proposition}

\section{Boundedness of ${\bold R}_r$.}\label{RBoundedness}

\indent
From the construction of operator ${\bold R}_r$ we conclude that in order
to prove necessary estimates for operator ${\bold R}_r$ it suffices
to prove these estimates for operator $R_r$. In the proposition below
we state necessary estimates for operator $R_r$.\\

\begin{proposition}\label{REstimate}
Let $0<\alpha< 1$, ${\bold M} \subset {\bold G}$ be a $C^{\infty}$ regular
q-concave CR submanifold of the form (\arabic{manifold}) and let $g \in \Gamma^{p, \alpha}_{(0,r)}({\bold M})$ be a form with compact support in
$U = {\widetilde U} \cap {\bold M}$.\\
\indent
Then operator $R_r$, defined in (\arabic{LocalFormula}) satisfies the
following estimate
$$\parallel R_r(g) \parallel_{\Gamma^{p,\alpha+1}_{(0,r-1)}(U)} <
C \cdot \parallel g \parallel_{\Gamma^{p, \alpha}_{(0,r)}(U)}
\eqno(\arabic{equation})
\newcounter{Restimate}
\setcounter{Restimate}{\value{equation}}
\addtocounter{equation}{1}$$
with a constant $C$ independent of $g$.
\end{proposition}

\indent
In our proof of proposition~\ref{REstimate} we will
use the approximation of $R_r$ by the operators
$$R_r(\epsilon)(f)(z) = (-1)^{r} \cdot \mbox{pr}_{\bold M} \circ
\frac{(n-1)!}{(2\pi i)^n}
\eqno(\arabic{equation})
\newcounter{OperatorREpsilon}
\setcounter{OperatorREpsilon}{\value{equation}}
\addtocounter{equation}{1}$$
$$\times \sum_{\iota}
\int_{{\bold M}_{\epsilon}\times [0,1]}
\vartheta_{\iota}(\zeta) \widetilde f(\zeta)
\wedge\omega^{\prime}_{r-1}\left((1-t)\frac{\bar\zeta - \bar z}
{{\mid \zeta - z \mid}^2} + t\frac{P^{\iota}(\zeta,z)}
{\Phi^{\iota}(\zeta,z)}\right) \wedge\omega(\zeta)$$
\noindent
when $\epsilon$ goes to $0$.\\
\indent
Using equality (\arabic{dbarFbar})
we obtain the following representation of kernels of these integrals
on ${\widetilde U} \times [0,1] \times {\bold M}:$
$$\left. \vartheta_{\iota}(\zeta) \cdot
\omega^{\prime}_{r-1}\left((1-t)\frac{\bar\zeta - \bar z}
{{\mid \zeta - z \mid}^2} + t\frac{P^{\iota}(\zeta,z)}
{\Phi^{\iota}(\zeta,z)}\right) \wedge\omega(\zeta)
\right|_{{\widetilde U} \times [0,1] \times {\bold M} }
\eqno(\arabic{equation})
\newcounter{cauchymartinellikernel}
\setcounter{cauchymartinellikernel}{\value{equation}}
\addtocounter{equation}{1}$$
$$= \sum_{i,J} a_{(i,J)}(t,\zeta,z)dt \wedge
\lambda^{i,J}_{r-1}(\zeta, z) +
\sum_{i,J} b_{(i,J)}(t,\zeta,z) dt \wedge
\gamma^{i,J}_{r-1}(\zeta, z),$$
\noindent
where $i$ is an index, $J= \cup_{i=1}^{11} J_i$ is a
multiindex such that $i \not \in J,$
$a_{(i,J)}(t,\zeta,z)$ and $b_{(i,J)}(t,\zeta,z)$ are polynomials in $t$
with coefficients that are smooth functions of  $z$, $\zeta$ and
$\theta(\zeta)$, and $\lambda^{i,J}_{r-1}(\zeta, z)$
and $\gamma^{i,J}_{r-1}(\zeta, z)$ are defined as follows:
$$\lambda^{i,J}_{r-1}(\zeta, z) =\frac{1}
{|\zeta - z|^{2(|J_1|+|J_7|+1)} \cdot
{\Phi(\zeta,z)}^{n-|J_1|-|J_7|-1}}
\eqno(\arabic{equation})
\newcounter{lambdaForm}
\setcounter{lambdaForm}{\value{equation}}
\addtocounter{equation}{1}$$
$$\times \sum \mbox{Det} \left[\bar{\zeta} - \bar{z},\hspace{0.03in}
Q^{(i)} {\overline F}^{(i)},\hspace{0.03in}
\overbrace{d\bar{\zeta}}^{j \in J_1},\hspace{0.03in}
\overbrace{Q d\rho}^{j \in J_2},\hspace{0.03in}
\overbrace{Q \cdot \chi}^{j \in J_3},\hspace{0.03in}
\overbrace{{\cal A}{\bar A} \cdot \bar\partial_{\zeta}a}^{j \in J_4},\right.$$
$$\left.\overbrace{a \cdot \mu_{\nu}}^{j \in J_5},\hspace{0.03in}
\overbrace{a \cdot \mu_{\tau}}^{j \in J_6}, \hspace{0.03in}
\overbrace{d\bar{z}}^{j \in J_7},\hspace{0.03in}
\overbrace{{\cal A}{\bar A} \cdot \bar \partial_z a}^{j \in J_8},
\hspace{0.03in}
\overbrace{a \cdot \bar\partial_z ({\cal A}{\bar A})}^{j \in J_9},
\hspace{0.03in}
\overbrace{{\overline F} \cdot \bar\partial_z Q}^{j \in J_{10}},
\hspace{0.03in}
\overbrace{Q \cdot \kappa}^{j \in J_{11}}
\right] \wedge\omega(\zeta),$$
\noindent
and
$$\gamma^{i,J}_{r-1}(\zeta, z) =\frac{1}
{|\zeta - z|^{2(|J_1|+|J_7|+1)} \cdot
{\Phi(\zeta,z)}^{n-|J_1|-|J_7|-1}}
\eqno(\arabic{equation})
\newcounter{gammaForm}
\setcounter{gammaForm}{\value{equation}}
\addtocounter{equation}{1}$$
$$\times \sum \mbox{Det} \left[\bar{\zeta} - \bar{z},\hspace{0.03in}
a_i {\cal A}{\bar A_i},\hspace{0.03in}
\overbrace{d\bar{\zeta}}^{j \in J_1},\hspace{0.03in}
\overbrace{Q d\rho}^{j \in J_2},\hspace{0.03in}
\overbrace{Q \cdot \chi}^{j \in J_3},\hspace{0.03in}
\overbrace{{\cal A}{\bar A} \cdot \bar\partial_{\zeta}a}^{j \in J_4},\right.$$
$$\left.\overbrace{a \cdot \mu_{\nu}}^{j \in J_5},\hspace{0.03in}
\overbrace{a \cdot \mu_{\tau}}^{j \in J_6}, \hspace{0.03in}
\overbrace{d\bar{z}}^{j \in J_7},\hspace{0.03in}
\overbrace{{\cal A}{\bar A} \cdot \bar \partial_z a}^{j \in J_8},
\hspace{0.03in}
\overbrace{a \cdot \bar\partial_z ({\cal A}{\bar A})}^{j \in J_9},
\hspace{0.03in}
\overbrace{{\overline F} \cdot \bar\partial_z Q}^{j \in J_{10}},
\hspace{0.03in}
\overbrace{Q \cdot \kappa}^{j \in J_{11}}
\right] \wedge\omega(\zeta).$$
\indent
In the proof of boundedness of operators $R_r$
we will need to know the differentiability  properties of integrals
with kernels $\lambda^{i,J}_{r-1}(\zeta, z)$ and
$\gamma^{i,J}_{r-1}(\zeta, z)$.
In the lemmas below we prepare necessary tools for the proof of
boundedness of $R_r$.\\
\indent
At first we consider the differentiability property of slightly different
kernels
$${\cal K}^{I}_{d,j}(\zeta, z) =
\frac{ \{ \rho(\zeta) \}^{I_1} (\zeta - z)^{I_2}
(\bar{\zeta} - \bar{z})^{I_3} }
{|\zeta -z|^d \cdot \Phi(\zeta,z)^{j}}
\overbrace{\wedge d_{\zeta}\mbox{Re}F^{(i)}}^{i \in I_4}
\overbrace{\wedge d\rho_i}^{i \in I_5}
\overbrace{\wedge d\theta_i(\zeta)}^{i \in I_6}
\wedge d\sigma_{2n-m}(\zeta),$$
\noindent
where $I= \cup_{j=1}^6 I_j$ and $I_j$ for $j=1, \dots, 6$ are
multiindices such that $I_1$ contains $m$ indices, $I_2$, $I_3$
contain $n$ indices, $I_4 \cup I_5 \cup I_6$ contains $m-1$ indices,
$|I_4|+|I_5|+|I_6|=m-1$, and $ \{ \rho(\zeta) \}^{I_1}=
\prod_{i_s \in I_1} \rho_s(\zeta)^{i_s}$,
$(\zeta - z)^{I_2} = \prod_{i_s \in I_2} (\zeta_s - z_s)^{i_s}$,
$(\bar{\zeta} - \bar{z})^{I_3} = \prod_{i_s \in I_3}
(\bar{\zeta}_s - \bar{z}_s)^{i_s}$.\\
\indent
For kernels ${\cal K}^{I}_{d,j}$ we introduce also the following notation
$$\begin{array}{llll}
k\left({\cal K}^{I}_{d,j}\right)=d-|I_2|-|I_3|,\vspace{0.1in}\\
h\left({\cal K}^{I}_{d,j}\right)=2j,\vspace{0.1in}\\
l\left({\cal K}^{I}_{d,j}\right)=|I_1|+|I_5|,\vspace{0.1in}\\
s\left({\cal K}^{I}_{d,j}\right)=|I_4|.
\end{array}$$

\begin{lemma}\label{Smoothness}
Let $U = {\widetilde U} \cap {\bold M}$ be a neighborhood
with $\Phi(\zeta, z)$ constructed in (\arabic{Barrier}),
$U(\epsilon) = {\widetilde U} \cap {\bold M}_{\epsilon}$ and let
$g(\zeta,z,\theta,t)$ be a smooth form with compact support in
${\widetilde U}_{\zeta} \times {\widetilde U}_{z} \times \S^{n-1}
\times[0,1]$.\\
\indent
Then for $g(\zeta, z, t)=g(\zeta,z,\theta(\zeta),t)$ and a vector field
$$D = \sum_{j=1}^n a_j(z) \frac{\partial}{\partial z_j} +
\sum_{j=1}^n b_j(z) \frac{\partial}{\partial\bar z_j}
\in {\bold C}T({{\bold M}})$$
the following equality holds
$$D \left( \int_{U(\epsilon) \times [0,1]} g(\zeta, z, t) \cdot
{\cal K}^{I}_{d,j}(\zeta, z) dt \right)
= \int_{U(\epsilon) \times [0,1]} \left[ D g(\zeta, z, t) \right]
\cdot {\cal K}^{I}_{d,j}(\zeta, z)dt
\eqno(\arabic{equation})
\newcounter{DzRepresentation}
\setcounter{DzRepresentation}{\value{equation}}
\addtocounter{equation}{1}$$
$$+ \int_{U(\epsilon) \times [0,1]}
\left[ D_{\zeta} g(\zeta, z, t) \right] \cdot
{\cal K}^{I}_{d,j}(\zeta, z) dt
+ \sum_{S,a,b} \int_{U(\epsilon) \times [0,1]}
c_{ \{ S,a,b \}}(\zeta, z, t) \cdot 
g(\zeta, z, t) \cdot {\cal K}^{S}_{a,b}(\zeta, z) dt$$
$$+ \sum_{k=1}^{m} \sum_{S} \int_{U(\epsilon) \times [0,1]}
c_{ \{ k, S \}}(\zeta, z, t) \cdot \left[ Y_{k,\zeta}(z) g(\zeta, z, t) \right]
\cdot {\cal K}^{S}_{d,j}(\zeta, z) dt$$
\noindent
where
$$g(\zeta,z,t)= g(\zeta,z,\theta(\zeta),t),$$
$$c_{ \{ S,a,b \}}(\zeta, z, t), c_{ \{ i, S \}}(\zeta, z, t)
\hspace{0.1in}\mbox{are}\hspace{0.1in}C^{\infty}
\hspace{0.1in}\mbox{functions of}\hspace{0.1in}\zeta,z,\theta(\zeta),t,$$
\noindent 
vector field $D_{\zeta}$ defined as
$$D_{\zeta} = \sum_{j=1}^n a_j(z) \frac{\partial}{\partial {\zeta}_j} +
\sum_{j=1}^n b_j(z) \frac{\partial}{\partial\bar{\zeta}_j}$$
and indices $a,b$ and multiindices $S$ satisfy the following conditions
$$\begin{array}{llll}
k\left({\cal K}^{S}_{a,b}\right)+h\left({\cal K}^{S}_{a,b}\right)
-l\left({\cal K}^{S}_{a,b}\right)-s\left({\cal K}^{S}_{a,b}\right)
\vspace{0.1in}\\
\leq k\left({\cal K}^{I}_{d,j}\right)+h\left({\cal K}^{I}_{d,j}\right)
-l\left({\cal K}^{I}_{d,j}\right)-s\left({\cal K}^{I}_{d,j}\right),
\vspace{0.1in}\\
k\left({\cal K}^{S}_{a,b}\right)+2h\left({\cal K}^{S}_{a,b}\right)
-2l\left({\cal K}^{S}_{a,b}\right)-2s\left({\cal K}^{S}_{a,b}\right)
\vspace{0.1in}\\
\leq k\left({\cal K}^{I}_{d,j}\right)+2h\left({\cal K}^{I}_{d,j}\right)
-2l\left({\cal K}^{I}_{d,j}\right)-2s\left({\cal K}^{I}_{d,j}\right).
\end{array}
\eqno(\arabic{equation})
\newcounter{SmoothnessIndices}
\setcounter{SmoothnessIndices}{\value{equation}}
\addtocounter{equation}{1}$$
\end{lemma}

\indent
{\bf Proof.}\\
\indent
To prove the lemma we represent integral from the left hand side of
(\arabic{DzRepresentation}) as
$$D \left( \int_{U(\epsilon) \times [0,1]} g(\zeta, z, t) \cdot
{\cal K}^{I}_{d,j}(\zeta, z) dt \right)
= \int_{U(\epsilon) \times [0,1]} \left[ D g(\zeta, z, t) \right]
\cdot {\cal K}^{I}_{d,j}(\zeta, z) dt
\eqno(\arabic{equation})
\newcounter{DzFormula}
\setcounter{DzFormula}{\value{equation}}
\addtocounter{equation}{1}$$
$$ - \int_{U(\epsilon) \times [0,1]} g(\zeta, z, t) 
\left[ D_{\zeta} {\cal K}^{I}_{d,j}(\zeta, z) \right]dt
+ \int_{U(\epsilon) \times [0,1]} g(\zeta, z, t) 
\left[ \left( D+ D_{\zeta} \right)
{\cal K}^{I}_{d,j}(\zeta, z) \right]dt.$$
\indent
To transform the second term of the right hand side of
(\arabic{DzFormula}) we apply integration by parts and obtain
$$\int_{U(\epsilon) \times [0,1]} g(\zeta, z, t) 
\left[ D_{\zeta} {\cal K}^{I}_{d,j}(\zeta, z) \right]dt
= - \int_{U(\epsilon) \times [0,1]} \left[ D_{\zeta} g(\zeta, z, t) \right]
{\cal K}^{I}_{d,j}(\zeta, z) dt.$$
\indent
To transform the third term of the right hand side of
(\arabic{DzFormula}) we will use the estimates below that follow from
the definitions of $F^{(k)}(\zeta,z)$ and ${\cal A}(\zeta,z)$ and
from the fact that $D \in {\bold C}T({{\bold M}})$
$$\begin{array}{lll}
\left( D + D_{\zeta} \right)
{\cal A}^2(\zeta, z) = {\cal O} \left( {\cal A}^2(\zeta, z) \right),
\vspace{0.2in}\\
\left( D + D_{\zeta} \right)
\mbox{Re} F^{(k)}(\zeta,z) = {\cal O} \left( |\zeta - z|^2 \right),
\vspace{0.2in}\\
\left( D + D_{\zeta} \right)
\mbox{Im} F^{(k)}(\zeta,z) = {\cal O} \left( |\zeta - z| \right).
\end{array}
\eqno(\arabic{equation})
\newcounter{DzDzetaFA}
\setcounter{DzDzetaFA}{\value{equation}}
\addtocounter{equation}{1}$$
Applying the operators $D$ and $D_{\zeta}$ to
${\cal K}^{I}_{d,j}(\zeta, z)$ we obtain
$$\left( D + D_{\zeta} \right) {\cal K}^{I}_{d,j}(\zeta, z)
\eqno(\arabic{equation})
\newcounter{DzDzetaK}
\setcounter{DzDzetaK}{\value{equation}}
\addtocounter{equation}{1}$$
$$= (-j) \frac{ \{ \rho(\zeta) \}^{I_1} (\zeta - z)^{I_2}
(\bar{\zeta} - \bar{z})^{I_3} }
{|\zeta -z|^d \cdot \Phi(\zeta,z)^{j+1}} \cdot
\left( \left[ \left( D + D_{\zeta} \right) {\cal A}^2 \right]
\overbrace{d_{\zeta}\mbox{Re}F^{(i)}}^{i \in I_4}
\wedge \overbrace{d\theta_i(\zeta)}^{i \in I_5}
\wedge d\sigma_{2n-m}(\zeta) \right.$$
$$+ \sum_{k=1}^m 2\mbox{Re}F^{(k)} \cdot
\left[ \left( D + D_{\zeta} \right) \mbox{Re}F^{(k)} \right]
\overbrace{d_{\zeta}\mbox{Re}F^{(i)}}^{i \in I_4}
\wedge \overbrace{d\theta_i(\zeta)}^{i \in I_5}
\wedge d\sigma_{2n-m}(\zeta)$$
$$\left. + \sum_{k=1}^m 2\mbox{Im}F^{(k)} \cdot
\left[ \left( D + D_{\zeta} \right) \mbox{Im}F^{(k)} \right]
\overbrace{d_{\zeta}\mbox{Re}F^{(i)}}^{i \in I_4}
\wedge \overbrace{d\theta_i(\zeta)}^{i \in I_5}
\wedge d\sigma_{2n-m}(\zeta) \right).$$
\indent
From estimates (\arabic{DzDzetaFA}) we conclude that the first two
terms of the right hand side of (\arabic{DzDzetaK})
can be represented as linear combinations with
$C^{\infty}\left( {\widetilde U}_{\zeta} \times {\widetilde U}_{z}
\times \S^{n-1} \times [0,1] \right)$
coefficients of kernels ${\cal K}^{S}_{d,j+1}$ with either
$$|S_2| + |S_3| = |I_2| + |I_3| + 4$$
or
$$|S_1| = |I_1| + 1 \hspace{0.1in}\mbox{and}\hspace{0.1in}
|S_2| + |S_3| = |I_2| + |I_3| + 2.$$
\indent
To handle the last term of the right hand side of (\arabic{DzDzetaK}) we
use the third estimate from (\arabic{DzDzetaFA}) and represent
corresponding integrals as
$$\int_{U(\epsilon) \times [0,1]}
(-j) \frac{ \{ \rho(\zeta) \}^{I_1} (\zeta - z)^{I_2}
(\bar{\zeta} - \bar{z})^{I_3}
2\mbox{Im}F^{(k)} }
{|\zeta -z|^d \cdot \Phi(\zeta,z)^{j+1}}
\eqno(\arabic{equation})
\newcounter{ImIntegral1}
\setcounter{ImIntegral1}{\value{equation}}
\addtocounter{equation}{1}$$
$$\times\left[ \left( D + D_{\zeta} \right) \mbox{Im}F^{(k)} \right]
\overbrace{d_{\zeta}\mbox{Re}F^{(i)}}^{i \in I_4}
\wedge \overbrace{d\theta_i(\zeta)}^{i \in I_5}
\wedge d\sigma_{2n-m}(\zeta)dt$$
$$= \sum_{ \{ |S_2|+|S_3|=|I_2|+|I_3|+1 \} } \int_{U(\epsilon) \times [0,1]}
c_{ \{ S \}}(\zeta, z, t) \cdot g(\zeta, z, t)
\frac{ \{ \rho(\zeta) \}^{I_1} (\zeta - z)^{S_2}
(\bar{\zeta} - \bar{z})^{S_3} } {|\zeta -z|^d}$$
$$\times d_{\zeta} \left(\frac{1}{\Phi(\zeta,z)}\right)^j
\overbrace{d_{\zeta}\mbox{Re}F^{(i)}}^{i \in I_4}
\wedge \overbrace{d\theta_i(\zeta)}^{i \in I_5}
\wedge \left( d_{\zeta}\mbox{Im}F^{(i)} \vee d\sigma_{2n-m}(\zeta)\right)dt$$
$$+ j \cdot \sum_{ \{ |S_2|+|S_3|=|I_2|+|I_3|+1 \} }
\int_{U(\epsilon) \times [0,1]}
c_{ \{ S \}}(\zeta, z, t) \cdot g(\zeta, z, t)
\frac{ \{ \rho(\zeta) \}^{I_1} (\zeta - z)^{S_2}
(\bar{\zeta} - \bar{z})^{S_3} } {|\zeta -z|^d}$$
$$\times \frac{ \sum_{i=k}^m 2 \mbox{Re}F^{(k)} d_{\zeta}\mbox{Re}F^{(k)} +
d_{\zeta}{\cal A}^2 } { \Phi(\zeta,z)^{j+1} }
\overbrace{d_{\zeta}\mbox{Re}F^{(i)}}^{i \in I_4}
\wedge \overbrace{d\theta_i(\zeta)}^{i \in I_5}
\wedge \left( d_{\zeta}\mbox{Im}F^{(k)} \vee d\sigma_{2n-m}(\zeta)\right)dt.$$
\indent
Kernels of the second term of the right hand side of (\arabic{ImIntegral1})
may be represented as linear combinations with
$C^{\infty}\left( {\widetilde U}_{\zeta} \times {\widetilde U}_{z}
\times \S^{n-1} \times [0,1] \right)$
coefficients of kernels ${\cal K}^{S}_{d,j+1}$ with  either
$$|S_2| + |S_3| = |I_2| + |I_3| + 4$$
or
$$|S_1| = |I_1| + 1 \hspace{0.1in}\mbox{and}\hspace{0.1in}
|S_2| + |S_3| = |I_2| + |I_3| + 2.$$
\indent
The first term of the right hand side of (\arabic{ImIntegral1}) we transform
applying integration by parts
$$\int_{U(\epsilon) \times [0,1]} c_{ \{ S \}}(\zeta, z, t)
\cdot g(\zeta, z, t) \frac{ \{ \rho(\zeta) \}^{I_1}(\zeta - z)^{S_2}
(\bar{\zeta} - \bar{z})^{S_3} }
{|\zeta -z|^d}
\eqno(\arabic{equation})
\newcounter{ImIntegral2}
\setcounter{ImIntegral2}{\value{equation}}
\addtocounter{equation}{1}$$
$$\times d_{\zeta} \left(\frac{1}{\Phi(\zeta,z)} \right)^j
\overbrace{d_{\zeta}\mbox{Re}F^{(i)}}^{i \in I_4}
\wedge \overbrace{d\theta_i(\zeta)}^{i \in I_5}
\wedge \left( d_{\zeta}\mbox{Im}F^{(k)} \vee d\sigma_{2n-m}(\zeta)\right)dt$$
$$= \int_{U(\epsilon) \times [0,1]}c_{ \{i,S \}}(\zeta, z, t) \frac{
\{ \rho(\zeta) \}^{I_1} (\zeta - z)^{S_2} (\bar{\zeta} - \bar{z})^{S_3} }
{|\zeta -z|^d \cdot \Phi(\zeta,z)^j }$$
$$\times \left[ Y_{k,\zeta}(z) g(\zeta, z, t) \right]
\overbrace{d_{\zeta}\mbox{Re}F^{(i)}}^{i \in I_4}
\wedge \overbrace{d\theta_i(\zeta)}^{i \in I_5}
\wedge d\sigma_{2n-m}(\zeta)dt$$
$$+ \sum_{ \{ d+1-|S_2|-|S_3| = d-|I_2|-|I_3| \} }
\int_{U(\epsilon) \times [0,1]}
c_{ \{ S \}}(\zeta, z, t) \cdot g(\zeta, z, t)
{\cal K}^{S}_{d+1,j}(\zeta, z) dt.$$
\indent
From (\arabic{DzDzetaK}), (\arabic{ImIntegral1}) and (\arabic{ImIntegral2})
follows statement of the lemma for the third term of the right
hand side of (\arabic{DzFormula}).\qed

\indent
In the next two lemmas we prove transformation and differentiation
formulas for the kernels
$${\cal B}^{T,I}_{d,j}(\zeta, z) =
\{\mbox{Im}F(\zeta, z) \}^{T} \cdot {\cal K}^{I}_{d,j}(\zeta, z),$$
and
$${\cal B}^{T,I}_{d,0}(\zeta, z) =
\{\mbox{Im}F(\zeta, z) \}^{T} \cdot {\cal K}^{I}_{d,0}(\zeta, z)
\cdot \log{\Phi(\zeta,z)},$$
where $\{\mbox{Im}F(\zeta, z) \}^{T} = \prod_{t_s \in T}
\left( \mbox{Im}F^{(s)}(\zeta, z) \right)^{t_s}$ for a multiindex
$T = (t_1, \dots, t_m)$.\\
\indent
Generalizing corresponding notations for ${\cal K}^{I}_{d,j}$ we denote
$$\begin{array}{llll}
k\left({\cal B}^{T,I}_{d,j}\right)=d-|I_2|-|I_3|,\vspace{0.1in}\\
h\left({\cal B}^{T,I}_{d,j}\right)=2j-|T|,\vspace{0.1in}\\
l\left({\cal B}^{T,I}_{d,j}\right)=|I_1|+|I_5|,\vspace{0.1in}\\
s\left({\cal B}^{T,I}_{d,j}\right)=|I_4|.
\end{array}$$

\begin{lemma}\label{Transformation}
Let $U = {\widetilde U} \cap {\bold M}$ be a neighborhood of
$z_0 \in {\bold M}$ with $\Phi(\zeta, z)$ constructed in (\arabic{Barrier}),
$U(\epsilon) = {\widetilde U} \cap {\bold M}_{\epsilon}$,
$g(\zeta,z,\theta,t)$ be a smooth form
with compact support and let $T$ be a
multiindex such that $t_k$ is a first nonzero index in $T$.\\
\indent
Then the following equality holds for
$g(\zeta, z, t)=g(\zeta,z,\theta(\zeta),t)$ and $j > 0$
$$\int_{U(\epsilon) \times [0,1]}
g(\zeta, z, t) \cdot {\cal B}^{T,I}_{d,j}(\zeta, z) dt
\eqno(\arabic{equation})
\newcounter{Transform}
\setcounter{Transform}{\value{equation}}
\addtocounter{equation}{1}$$
$$= c \cdot \int_{U(\epsilon) \times [0,1]}
\left[ Y_{k,\zeta}(z) g(\zeta, z, t) \right] \cdot
{\cal B}^{\hat{T},I}_{d,j-1}(\zeta, z) dt$$
$$+ \sum_{ \left\{ L, S, a, b \right\} }
\int_{U(\epsilon) \times [0,1]} c_{\left\{ L, S, a, b \right\} }(\zeta, z, t)
\cdot g(\zeta, z, t) \cdot {\cal B}^{L,S}_{a,b}(\zeta, z) dt,$$
\noindent
where $\hat{T} = (t_1, \cdots, t_k -1, \cdots, t_m)$, $c$ is a constant,
$c_{ \{ L, S, a, b \}}(\zeta,z,t)
=c_{ \{ L, S, a, b \}}(\zeta,z,\theta(\zeta),t)$ with
$$c_{ \{ L, S, a, b \}}(\zeta,z,\theta,t)
\in C^{\infty}\left( {\widetilde U}_{\zeta} \times {\widetilde U}_{z}
\times \S^{n-1} \times [0,1] \right)$$
\noindent
and indices $a, b$ and multiindices $L, S$ are such that
$$\begin{array}{lllll}
|L| < |T|,\vspace{0.1in}\\
k\left({\cal B}^{L,S}_{a,b}\right)+h\left({\cal B}^{L,S}_{a,b}\right)
-l\left({\cal B}^{L,S}_{a,b}\right)-s\left({\cal B}^{L,S}_{a,b}\right)
\vspace{0.1in}\\
\leq k\left({\cal B}^{T,I}_{d,j}\right)+h\left({\cal B}^{T,I}_{d,j}\right)
-l\left({\cal B}^{T,I}_{d,j}\right)-s\left({\cal B}^{T,I}_{d,j}\right),
\vspace{0.1in}\\
k\left({\cal B}^{L,S}_{a,b}\right)+2h\left({\cal B}^{L,S}_{a,b}\right)
-2l\left({\cal B}^{L,S}_{a,b}\right)-2s\left({\cal B}^{L,S}_{a,b}\right)
\vspace{0.1in}\\
\leq k\left({\cal B}^{T,I}_{d,j}\right)+2h\left({\cal B}^{T,I}_{d,j}\right)
-2l\left({\cal B}^{T,I}_{d,j}\right)-2s\left({\cal B}^{T,I}_{d,j}\right).
\end{array}
\eqno(\arabic{equation})
\newcounter{TransformIndices}
\setcounter{TransformIndices}{\value{equation}}
\addtocounter{equation}{1}$$
\indent
For $j=1$ and $|T|=1$ stronger inequalities hold for the terms
${\cal B}^{\emptyset,S}_{a,0}$ of the right hand side of (\arabic{Transform})
$$\begin{array}{llll}
k\left({\cal B}^{\emptyset,S}_{a,0}\right)
-l\left({\cal B}^{\emptyset,S}_{a,0}\right)
-s\left({\cal B}^{\emptyset,S}_{a,0}\right)+1
\vspace{0.1in}\\
\leq k\left({\cal B}^{T,I}_{d,1}\right)+h\left({\cal B}^{T,I}_{d,1}\right)
-l\left({\cal B}^{T,I}_{d,1}\right)-s\left({\cal B}^{T,I}_{d,1}\right),
\vspace{0.1in}\\
k\left({\cal B}^{\emptyset,S}_{a,0}\right)
-2l\left({\cal B}^{\emptyset,S}_{a,0}\right)
-2s\left({\cal B}^{\emptyset,S}_{a,0}\right)+2
\vspace{0.1in}\\
\leq k\left({\cal B}^{T,I}_{d,1}\right)+2h\left({\cal B}^{T,I}_{d,1}\right)
-2l\left({\cal B}^{T,I}_{d,1}\right)-2s\left({\cal B}^{T,I}_{d,1}\right).
\end{array}
\eqno(\arabic{equation})
\newcounter{LogTransformIndices}
\setcounter{LogTransformIndices}{\value{equation}}
\addtocounter{equation}{1}$$
\end{lemma}

\indent
{\bf Proof.}\\
\indent
The integral in the left hand side of (\arabic{Transform}) may be
represented as a sum of two integrals
$$\int_{U(\epsilon) \times [0,1]} g(\zeta, z, t) \cdot
\{\mbox{Im}F(\zeta, z) \}^{T}{\cal K}^{I}_{d,j}(\zeta, z) dt
\eqno(\arabic{equation})
\newcounter{Transform1}
\setcounter{Transform1}{\value{equation}}
\addtocounter{equation}{1}$$
$$= - \frac{1}{2(j-1)} \int_{U(\epsilon) \times [0,1]} g(\zeta, z, t)
\{\mbox{Im}F(\zeta, z) \}^{\hat{T}}
\frac{ \{ \rho(\zeta) \}^{I_1}
(\zeta - z)^{I_2} (\bar{\zeta} - \bar{z})^{I_3} } {|\zeta -z|^d}$$
$$\times d_{\zeta} \left(\frac{1}{\Phi(\zeta,z)}\right)^{j-1}
\overbrace{d_{\zeta}\mbox{Re}F^{(i)}}^{i \in I_4}
\wedge \overbrace{d\theta_i(\zeta)}^{i \in I_5}
\wedge \left( d_{\zeta}\mbox{Im}F^{(k)} \vee d\sigma_{2n-m}(\zeta)\right)dt$$
$$ - \int_{U(\epsilon) \times [0,1]} g(\zeta, z, t)
\{\mbox{Im}F(\zeta, z) \}^{\hat{T}} \frac{ \{ \rho(\zeta) \}^{I_1}
(\zeta - z)^{I_2} (\bar{\zeta} - \bar{z})^{I_3} }
{|\zeta -z|^d}$$
$$\times \frac{ \sum_{i=1}^m \mbox{Re}F^{(i)} d_{\zeta}\mbox{Re}F^{(i)} +
\frac{1}{2} d_{\zeta}{\cal A}^2 }{\Phi(\zeta,z)^{j}}
\overbrace{d_{\zeta}\mbox{Re}F^{(i)}}^{i \in I_4}
\wedge \overbrace{d\theta_i(\zeta)}^{i \in I_5}
\wedge \left( d_{\zeta}\mbox{Im}F^{(k)} \vee d\sigma_{2n-m}(\zeta)\right)dt.$$
\indent
As follows from (\arabic{DzDzetaFA}) the second integral admits necessary
representation with indices $a,b$ and multiindices $L, S$ satisfying
conditions (\arabic{TransformIndices}).
To transform the first integral from the right hand side of
(\arabic{Transform1}) we apply integration by parts and obtain
$$- \frac{1}{2(j-1)} \int_{U(\epsilon) \times [0,1]} g(\zeta, z, t)
\{\mbox{Im}F(\zeta, z) \}^{\hat{T}}
\frac{ \{ \rho(\zeta) \}^{I_1} (\zeta - z)^{I_2}
(\bar{\zeta} - \bar{z})^{I_3} } {|\zeta -z|^d}$$
$$\times d_{\zeta} \left(\frac{1}{\Phi(\zeta,z)}\right)^{j-1}
\overbrace{d_{\zeta}\mbox{Re}F^{(i)}}^{i \in I_4}
\wedge \overbrace{d\theta_i(\zeta)}^{i \in I_5}
\wedge \left( d_{\zeta}\mbox{Im}F^{(k)} \vee d\sigma_{2n-m}(\zeta)\right)dt$$
$$= \frac{1}{2(j-1)} \left[ \int_{U(\epsilon) \times [0,1]}
\left[ Y_{k,\zeta}(z) g(\zeta, z, t) \right] \cdot
\{\mbox{Im}F(\zeta, z) \}^{\hat{T}}{\cal K}^{I}_{d,j-1}(\zeta, z) dt \right.$$
$$+ \sum_{ \left\{|L| = |T|-2 \right\} }
\int_{U(\epsilon) \times [0,1]}
c_{\left\{ L \right\}}(\zeta, z, t)\cdot g(\zeta, z, t) \cdot
\{\mbox{Im}F(\zeta, z) \}^{L} {\cal K}^{I}_{d,j-1}(\zeta, z) dt$$
$$\left.+ \sum_{ \left\{|S_3|+|S_4|=|I_2|+|I_3|+1 \right\} }
\int_{U(\epsilon) \times [0,1]}
c_{ \left\{ S \right\}}(\zeta, z, t) \cdot
g(\zeta, z, t) \{\mbox{Im}F(\zeta, z) \}^{\hat{T}}
{\cal K}^{S}_{d+2,j-1}(\zeta, z) dt\right]$$
with kernels satisfying (\arabic{TransformIndices}).\qed

\begin{lemma}\label{LogSmoothness}
Let $U = {\widetilde U} \cap {\bold M}$ be a neighborhood of
$z_0 \in {\bold M}$ with $\Phi(\zeta, z)$ constructed in (\arabic{Barrier}),
$U(\epsilon) = {\widetilde U} \cap {\bold M}_{\epsilon}$, and let
$g(\zeta,z,\theta,t)$ be a smooth form with compact support in
${\widetilde U}_{\zeta} \times {\widetilde U}_{z} \times \S^{n-1}
\times[0,1]$.\\
\indent
Then for $g(\zeta, z, t)=g(\zeta,z,\theta(\zeta),t)$ and a vector field
$$D = \sum_{j=1}^n a_j(z) \frac{\partial}{\partial z_j} +
\sum_{j=1}^n b_j(z) \frac{\partial}{\partial\bar z_j}
\in {\bold C}T({{\bold M}})$$
the following equality holds
$$D \int_{U(\epsilon) \times [0,1]}
g(\zeta, z, t) \cdot {\cal B}^{\emptyset,I}_{d,0}(\zeta, z) dt
= \int_{U(\epsilon) \times [0,1]}
Dg(\zeta, z, t) \cdot {\cal B}^{\emptyset,I}_{d,0}(\zeta, z) dt
\eqno(\arabic{equation})
\newcounter{DLog}
\setcounter{DLog}{\value{equation}}
\addtocounter{equation}{1}$$
$$+ \int_{U(\epsilon) \times [0,1]} \left[ D_{\zeta} g(\zeta, z, t) \right]
\cdot {\cal B}^{\emptyset,I}_{d,0}(\zeta, z) dt
+ \sum_{ \left\{ S, a \right\} }
\int_{U(\epsilon) \times [0,1]} c_{\left\{ S, a \right\} }(\zeta, z, t)
g(\zeta, z, t) \cdot {\cal K}^{S}_{a,1}(\zeta, z) dt$$
$$+ \sum_{k=1}^m \sum_{ \left\{ S \right\} }
\int_{U(\epsilon) \times [0,1]} c_{\left\{ S \right\} }(\zeta, z, t)
\left[ Y_{k,\zeta}(z) g(\zeta, z, t) \right] \cdot
{\cal B}^{\emptyset,S}_{d,0}(\zeta, z) dt,$$
with $a$ and $S$ satisfying (\arabic{SmoothnessIndices}).
\end{lemma}

\indent
{\bf Proof.}\\
\indent
Proof of the lemma is analogous to the proof of lemma~\ref{Smoothness}.\\
\indent
We represent the integral from the left hand side of (\arabic{DLog}) as
$$D \left( \int_{U(\epsilon) \times [0,1]} g(\zeta, z, t) \cdot
{\cal B}^{\emptyset,I}_{d,0}(\zeta, z) dt \right)
= \int_{U(\epsilon) \times [0,1]} \left[ D g(\zeta, z, t) \right]
\cdot {\cal B}^{\emptyset,I}_{d,0}(\zeta, z) dt
\eqno(\arabic{equation})
\newcounter{DLogFormula}
\setcounter{DLogFormula}{\value{equation}}
\addtocounter{equation}{1}$$
$$ - \int_{U(\epsilon) \times [0,1]} g(\zeta, z, t) 
\left[ D_{\zeta} {\cal B}^{\emptyset,I}_{d,0}(\zeta, z) \right]dt
+ \int_{U(\epsilon) \times [0,1]} g(\zeta, z, t) 
\left[ \left( D+ D_{\zeta} \right)
{\cal B}^{\emptyset,I}_{d,0}(\zeta, z) \right]dt.$$
\indent
To transform the second term of the right hand side of
(\arabic{DLogFormula}) we apply integration by parts and obtain
$$\int_{U(\epsilon) \times [0,1]} g(\zeta, z, t) 
\left[ D_{\zeta} {\cal B}^{\emptyset,I}_{d,0}(\zeta, z) \right]dt
= - \int_{U(\epsilon) \times [0,1]} \left[ D_{\zeta} g(\zeta, z, t) \right]
{\cal B}^{\emptyset,I}_{d,0}(\zeta, z) dt.$$
\indent
To transform the third term of the right hand side of
(\arabic{DLogFormula}) we use estimates (\arabic{DzDzetaFA}) and obtain
$$\left( D + D_{\zeta} \right) {\cal B}^{\emptyset,I}_{d,0}(\zeta, z)
= \sum_{ \left\{ S, a \right\} }
c_{\left\{ S, a \right\} }(\zeta, z, t)
\cdot {\cal K}^{S}_{a,1}(\zeta, z)
\eqno(\arabic{equation})
\newcounter{DzDzetaL}
\setcounter{DzDzetaL}{\value{equation}}
\addtocounter{equation}{1}$$
$$+ \sum_{k=1}^m \frac{ \{ \rho(\zeta) \}^{I_1} (\zeta - z)^{I_2}
(\bar{\zeta} - \bar{z})^{I_3} }
{|\zeta -z|^d \cdot \Phi(\zeta,z)}$$
$$\times 2\mbox{Im}F^{(k)} \cdot
\left[ \left( D + D_{\zeta} \right) \mbox{Im}F^{(k)} \right]
\overbrace{d_{\zeta}\mbox{Re}F^{(i)}}^{i \in I_4}
\wedge \overbrace{d\theta_i(\zeta)}^{i \in I_5}
\wedge d\sigma_{2n-m}(\zeta).$$
\indent
Again using estimates (\arabic{DzDzetaFA}) and applying integration
by parts as in lemma~\ref{Smoothness} we obtain (\arabic{DLog}).\qed

\indent
The following two simple lemmas (cf. \cite{P}) will be used in the
further estimates.\\

\begin{lemma}\label{TangentDifferentials}
Let ${\bold M}$ be a generic CR submanifold in the unit ball
${\bold B}^n$ in $\C^n$ of the form:
$${\bold M} = \{ z \in {\bold B}^n: \rho_{1}(z) = \dots = \rho_{m}(z) = 
0\},$$
\noindent
where $\{ \rho_{k} \}, \  k = 1, \dots , m  \ ( m < n)$ are real
valued functions of the class $C^{\infty}$ satisfying
$$\partial \rho_1 \wedge \cdots \wedge \partial \rho_m \neq 0 
\hspace{0.3in} \mbox{on}{\bold M}.$$
\indent
Then for any point $\zeta_0 \in {\bold M}$ there exists a neighborhood
${\bold V}_{\epsilon}(\zeta_0) = \{ \zeta : | \zeta - \zeta_0 | <
\epsilon \}$ such that for any $n \geq s>n-m$ and $p>2n-s-m$ the
following representation holds in ${\bold V}_{\epsilon}:$
$$d{\bar\zeta}_{i_1} \wedge \dots d{\bar\zeta}_{i_p} \wedge
d \zeta_{k_1} \wedge \dots d \zeta_{k_s} = \sum d \rho_{j_1} \wedge
\dots d \rho_{j_{p-(2n-s-m)}} \wedge
g_{j_1 \dots j_{p-(2n-s-m)}}^{i_1 \dots i_p}(\zeta)
\eqno(\arabic{equation})
\newcounter{tangentestimate}
\setcounter{tangentestimate}{\value{equation}}
\addtocounter{equation}{1}$$
with $g_{j_1 \dots j_{p-(2n-s-m)}}^{i_1 \dots i_p}$ of the class
$C^{\infty}({\bold V}_{\epsilon})$.\qed
\end{lemma}

\begin{lemma}\label{Integral}
Let
$${\bold B}(1) = \{ (\rho, \eta, w) \in \R^s \times \R^m \times
\C^{n-m}: \sum_{i=1}^s \rho_i^2 + \sum_{i=1}^m \eta_i^2
+ \sum_{i=1}^{n-m} |w|^2 < 1 \},$$
$${\bold V}(\delta) = \{ (\rho, \eta, w) \in \R^s \times \R^m \times
\C^{n-m}: \sum_{i=1}^s |\rho_i| + \sum_{i=1}^m |\eta_i|
+ \sum_{i=1}^{n-m} |w|^2 < \delta^2 \},$$
$$K\left\{ \alpha, k, h, s\right\}(\rho, \eta, w, \epsilon)
= \frac{1}{ (\epsilon + \sum_{i=1}^s |\rho_i|+\sum_{i=1}^{m}|\eta_i|
+ \sum_{i=1}^{n-m} |w_i|)^{k} }$$
$$\times \frac{ \wedge_{i=1}^{s} d\rho_i \wedge_{i=1}^{m} d\eta_i
\wedge_{i=1}^{n-m} (dw_i \wedge d\bar{w}_i)}
{(\sqrt{\epsilon} + \sum_{i=1}^s \sqrt{|\rho_i|}
+ \sum_{i=1}^m \sqrt{|\eta_i|} + \sum_{i=1}^{n-m}
|w_i|)^{2h-\alpha}},$$
\noindent
with $0 \leq \alpha <1$ and $k, h, s \in \Z.$\\
\indent
Let
$${\cal I}_1\left\{\alpha,k,h,s\right\}(\epsilon,\delta)
= \int_{{\bold V}(\delta)} K\left\{ \alpha,k,h,s\right\}(\eta,w,\epsilon),$$
\noindent
and
$${\cal I}_2\left\{\alpha,k,h,s\right\}(\epsilon,\delta)
= \int_{{\bold B}(1) \setminus {\bold V}(\delta)}
K\left\{ \alpha,k,h,s\right\}(\eta, w, \epsilon).$$

\indent
Then
$${\cal I}_1\left\{0,k,h,s\right\}(\epsilon,\delta)$$
$$= \left\{ \begin{array}{ll}
\left\{ \begin{array}{ll}
{\cal O}\left(\epsilon^{2n-m-k-h+s} \cdot (\log{\epsilon})^2 \right)
& \mbox{if} \hspace{0.05in} k \geq 2n-2m \hspace{0.05in}
\mbox{and} \hspace{0.05in} k+h-s \geq 2n-m,\\
{\cal O}\left(\delta \right)
& \mbox{if} \hspace{0.05in} k \geq 2n-2m \hspace{0.05in}
\mbox{and} \hspace{0.05in} k+h-s \leq 2n-m-1,\\
\end{array} \right. \vspace{0.1in}\\
\left\{ \begin{array}{ll}
{\cal O}\left( \epsilon^{(2n-k-2(h-s))/2} \cdot \log{\epsilon}\right) &
\mbox{if} \hspace{0.05in}
k \leq 2n-2m-1 \hspace{0.05in} \mbox{and} \hspace{0.05in} k+2(h-s) \geq 2n,\\
{\cal O}\left( \delta \right) & \mbox{if} \hspace{0.05in}
k \leq 2n-2m-1 \hspace{0.05in} \mbox{and} \hspace{0.05in} k+2(h-s) \leq 2n-1,
\end{array} \right.
\end{array} \right.$$

$${\cal I}_1\left\{\alpha,k,h,s\right\}(\epsilon, \delta)$$
$$= {\cal O}\left(\delta^{\alpha} \right)\hspace{0.1in}\mbox{if}
\hspace{0.1in} \alpha > 0,\hspace{0.1in}
\mbox{and}
\hspace{0.1in}\left\{
\begin{array}{ll}
\mbox{if} \hspace{0.05in}
k \geq 2n-2m \hspace{0.05in}\mbox{and}\hspace{0.05in} k+h-s \leq 2n-m-1\\
\mbox{if} \hspace{0.05in}
k \leq 2n-2m-1 \hspace{0.05in}\mbox{and}\hspace{0.05in} k+2(h-s) \leq 2n,
\end{array} \right.$$

$${\cal I}_2\left\{\alpha,k,h,s\right\}(\epsilon, \delta)$$
$$= {\cal O}\left(\delta^{\alpha-1} \right)\hspace{0.1in}\mbox{if}
\hspace{0.1in} \alpha > 0,\hspace{0.1in}
\mbox{and}
\hspace{0.1in}\left\{
\begin{array}{ll}
\mbox{if} \hspace{0.05in}
k \geq 2n-2m \hspace{0.05in}\mbox{and}\hspace{0.05in} k+h-s \leq 2n-m\\
\mbox{if} \hspace{0.05in}
k \leq 2n-2m-1 \hspace{0.05in}\mbox{and}\hspace{0.05in} k+2(h-s) \leq 2n+1,
\end{array} \right.$$

\noindent
and
$${\cal I}_2\left\{\alpha,k,h,s\right\}(\epsilon, \delta)$$
$$= {\cal O}\left(\delta^{\alpha-2} \right)\hspace{0.1in}\mbox{if}
\hspace{0.1in} \alpha > 0,\hspace{0.1in}
\mbox{and}
\hspace{0.1in}\left\{
\begin{array}{ll}
\mbox{if} \hspace{0.05in}
k \geq 2n-2m \hspace{0.05in}\mbox{and}\hspace{0.05in} k+h-s \leq 2n-m\\
\mbox{if} \hspace{0.05in}
k \leq 2n-2m-1 \hspace{0.05in}\mbox{and}\hspace{0.05in} k+2(h-s) \leq 2n+2.
\end{array} \right.$$
\qed
\end{lemma}

\indent
{\bf Proof of proposition~\ref{REstimate}.}\\

\indent
According to (\arabic{OperatorREpsilon}) in order to prove statement
of the proposition it suffices to prove the estimates
$$\begin{array}{ll}
\| \int_{U(\epsilon) \times [0,1]}
a_{(i,J)}(t,\zeta,z)dt \wedge {\widetilde g}(\zeta) \wedge
\lambda^{i,J}_{r-1}(\zeta, z) \|_{ \Gamma^{p,\alpha +1}(U) } \leq
C \cdot \| g \|_{ \Gamma^{p, \alpha}(U) },
\vspace{0.2in}\\
\| \int_{U(\epsilon) \times [0,1]}
b_{(i,J)}(t,\zeta,z)dt \wedge {\widetilde g}(\zeta) \wedge
\gamma^{i,J}_{r-1}(\zeta, z) \|_{ \Gamma^{p,\alpha +1}(U) } \leq
C \cdot \| g \|_{ \Gamma^{p, \alpha}(U) }
\end{array}
\eqno(\arabic{equation})
\newcounter{lambdagammaEstimates}
\setcounter{lambdagammaEstimates}{\value{equation}}
\addtocounter{equation}{1}$$
\noindent
with constant $C$ independent of $g$ and $\epsilon$.\\
\indent
Using the estimates
$$\begin{array}{llll}
| {\cal A}{\bar A} \cdot \bar\partial_{\zeta} a | = {\cal O}
(|\zeta - z|^3), \hspace{0.1in} |\mu_{\nu}| = {\cal O} (|\zeta - z|^3),
\hspace{0.1in} |\mu_{\tau}| = {\cal O} (|\zeta - z|^2),\vspace{0.1in}\\
| {\cal A}{\bar A} \cdot \bar\partial_z a | = {\cal O}
(|\zeta - z|^3), \hspace{0.1in}
| a \cdot \bar\partial_{z}({\cal  A}{\bar A}) | = {\cal O}(|\zeta - z|^2),
\vspace{0.1in}\\
|\kappa(\zeta,z)| = {\cal O}(|\zeta - z|),\hspace{0.1in}
|\chi(\zeta,z)| = {\cal O}(|\zeta - z|),
\vspace{0.1in}\\
F^{(k)}(\zeta,z) = \frac{1}{2} \left( \rho_k(z) - \rho_k(\zeta) \right)
+ {\cal O}(|\zeta - z|^2) + \sqrt{-1} \mbox{Im} F^{(k)}(\zeta,z)
\end{array}
\eqno(\arabic{equation})
\newcounter{termsestimates}
\setcounter{termsestimates}{\value{equation}}
\addtocounter{equation}{1}$$
for the terms of the determinants in
(\arabic{lambdaForm}) and (\arabic{gammaForm})
and applying lemma~\ref{TangentDifferentials} to the differential form
$$\overbrace{d \bar{\zeta}}^{|J_1|+r} \wedge
\overbrace{\chi}^{|J_3|} \wedge
\overbrace{\mu_{\tau}}^{|J_6|} \wedge
\omega(\zeta)$$
we obtain representations
$$a_{(i,J)}(t,\zeta,z) dt \wedge {\widetilde g}(\zeta) \wedge
\lambda^{i,J}_{r-1}(\zeta, z)$$
$$= \sum_{|T|+|E| \leq |J_{10}|+1}
c_{ \{I,d,j \} }(\zeta, z, t) {\widetilde g}(\zeta)
\{\mbox{Im}F(\zeta, z) \}^{T} {\cal K}^{I(J,S,T,E)}_{d,j}(\zeta, z),
\eqno(\arabic{equation})
\newcounter{lambdaRepresentation}
\setcounter{lambdaRepresentation}{\value{equation}}
\addtocounter{equation}{1}$$
and
$$b_{(i,J)}(t,\zeta,z) dt \wedge {\widetilde g}(\zeta) \wedge
\gamma^{i,J}_{r-1}(\zeta, z)$$
$$= \sum_{|T|+|E| \leq |J_{10}|}
c_{ \{I,d,j \} }(\zeta, z, t) {\widetilde g}(\zeta)
\{\mbox{Im}F(\zeta, z) \}^{T} {\cal K}^{I(J,S,T,E)}_{d,j}(\zeta, z).
\eqno(\arabic{equation})
\newcounter{gammaRepresentation}
\setcounter{gammaRepresentation}{\value{equation}}
\addtocounter{equation}{1}$$
\indent
Multiindices $T$ and $E$ in (\arabic{lambdaRepresentation}) are
obtained from the decomposition
$$\{ \overline{F}(\zeta, z) \}^{ \{J_{10} \cup i \} }$$
$$= \sum_{|T|+|E|+\frac{1}{2}\left(|G|+|H|\right)=|J_{10}|+1}
c_{\left\{ T,E,G,H \right\} }(\zeta, z)
\{ \mbox{Im}F(\zeta,z)\}^{T} \{ \rho(\zeta) \}^{E}
(\zeta - z)^{G} (\bar{\zeta} - \bar{z})^{H}$$
and multiindices $I_i$ for $i=1, \dots, 6$ and indices $d,j$ in
(\arabic{lambdaRepresentation}) satisfy the conditions below
$$\begin{array}{lllllll}
\hspace{0.3in}d = 2(|J_1|+|J_7|+1),\vspace{0.1in}\\
\hspace{0.3in}j = n-|J_1|-|J_7|-1,\vspace{0.1in}\\
\hspace{0.3in}|I_1| = |E|,\vspace{0.1in}\\
\hspace{0.3in}|I_2| + |I_3| =
1+|J_3|+3|J_4|+3|J_5|+2|J_6|+3|J_8|+2|J_9|+|J_{11}|
\vspace{0.1in}\\
+2\left( |J_{10}|+1-|T|-|E| \right),\vspace{0.1in}\\
\hspace{0.3in}|I_4| = 0,\vspace{0.1in}\\
\hspace{0.3in}|I_5| = |J_1|+|J_2|+|J_3|+|J_6|+r+m-n.
\end{array}
\eqno(\arabic{equation})
\newcounter{KlambdaIndices}
\setcounter{KlambdaIndices}{\value{equation}}
\addtocounter{equation}{1}$$
\indent
Multiindices $T$ and $E$ in (\arabic{gammaRepresentation}) are
obtained from the decomposition
$$\{ \overline{F}(\zeta, z) \}^{ \{J_{10} \} }$$
$$= \sum_{|T|+|E|+\frac{1}{2}\left(|G|+|H|\right)=|J_{10}|}
c_{\left\{ T,E,G,H \right\} }(\zeta, z)
\{ \mbox{Im}F(\zeta,z)\}^{T} \{ \rho(\zeta) \}^{E}
(\zeta - z)^{G} (\bar{\zeta} - \bar{z})^{H}$$
and multiindices $I_i$ for $i=1, \dots, 6$ and indices $d,j$ in
(\arabic{gammaRepresentation}) satisfy the conditions
$$\begin{array}{lllllll}
\hspace{0.3in}d = 2(|J_1|+|J_7|+1),\vspace{0.1in}\\
\hspace{0.3in}j = n-|J_1|-|J_7|-1,\vspace{0.1in}\\
\hspace{0.3in}|I_1| = |E|,\vspace{0.1in}\\
\hspace{0.3in}|I_2| + |I_3| =
4+|J_3|+3|J_4|+3|J_5|+2|J_6|+3|J_8|+2|J_9|+|J_{11}|
\vspace{0.1in}\\
+2\left( |J_{10}|-|T|-|E| \right),\vspace{0.1in}\\
\hspace{0.3in}|I_4| = 0,\vspace{0.1in}\\
\hspace{0.3in}|I_5| = |J_1|+|J_2|+|J_3|+|J_6|+r+m-n.
\end{array}
\eqno(\arabic{equation})
\newcounter{KgammaIndices}
\setcounter{KgammaIndices}{\value{equation}}
\addtocounter{equation}{1}$$
\indent
Using representations (\arabic{lambdaRepresentation}) and
(\arabic{gammaRepresentation}) we reduce the problem
of proving (\arabic{lambdagammaEstimates}) to each term
$${\widetilde g}(\zeta){\cal B}^{T,I}_{d,j}(\zeta, z)$$
of the right hand side of these representations.\\
\indent
We further reduce the problem using transformation from
lemma~\ref{Transformation} to obtain a representation
$$\int_{U(\epsilon) \times [0,1]}{\widetilde g}(\zeta)
{\cal B}^{T,I}_{d,j}(\zeta, z) dt
\eqno(\arabic{equation})
\newcounter{BRepresentation1}
\setcounter{BRepresentation1}{\value{equation}}
\addtocounter{equation}{1}$$
$$= \sum_{i=0}^{p} \sum_{|H|=i} \int_{U(\epsilon) \times [0,1]}
c_{\left\{ H, M \right\} }(\zeta, z, t)
\left[ \{ Y_{\zeta}(z) \}^{H} {\widetilde g}(\zeta) \right]
{\cal B}^{P(H),M}_{a, b}(\zeta, z) dt,$$
\noindent
where
$$\{ Y_{\zeta}(z) \}^{H} {\widetilde g}(\zeta)
:= Y_{h_1,\zeta}(z) \circ \cdots \circ Y_{h_p,\zeta}(z){\widetilde g}(\zeta),$$
\noindent
$c_{ \{ H, M \}}(\zeta,z,t) = c_{ \{ H, M \}}(\zeta,z,\theta(\zeta),t)$ with
$$c_{ \{ H, M \}}(\zeta,z,\theta,t)
\in C^{\infty}\left( {\widetilde U}_{\zeta} \times {\widetilde U}_{z}
\times \S^{n-1} \times [0,1] \right)$$
\noindent
and indices $a, b$ and multiindices $H, P(H), M$ are such that
$$\begin{array}{lllll}
P(H) = \emptyset \hspace{0.1in}\mbox{if}\hspace{0.1in}|H| < p,\vspace{0.1in}\\
k\left({\cal B}^{P(H),M}_{a,b}\right)+h\left({\cal B}^{P(H),M}_{a,b}\right)
-l\left({\cal B}^{P(H),M}_{a,b}\right)-s\left({\cal B}^{P(H),M}_{a,b}\right)
\vspace{0.05in}\\
\leq k\left({\cal B}^{T,I}_{d,j}\right)+h\left({\cal B}^{T,I}_{d,j}\right)
-l\left({\cal B}^{T,I}_{d,j}\right)-s\left({\cal B}^{T,I}_{d,j}\right)-|H|,
\vspace{0.1in}\\
k\left({\cal B}^{P(H),M}_{a,b}\right)+2h\left({\cal B}^{P(H),M}_{a,b}\right)
-2l\left({\cal B}^{P(H),M}_{a,b}\right)-2s\left({\cal B}^{P(H),M}_{a,b}\right)
\vspace{0.05in}\\
\leq k\left({\cal B}^{T,I}_{d,j}\right)+2h\left({\cal B}^{T,I}_{d,j}\right)
-2l\left({\cal B}^{T,I}_{d,j}\right)-2s\left({\cal B}^{T,I}_{d,j}\right)-2|H|.
\end{array}
\eqno(\arabic{equation})
\newcounter{1StepIndices}
\setcounter{1StepIndices}{\value{equation}}
\addtocounter{equation}{1}$$
\indent
To obtain such a representation we repeatedly apply transformation from
lemma~\ref{Transformation} to a selected term of the right hand side
of (\arabic{lambdaRepresentation}) or (\arabic{gammaRepresentation})
and all the resulting terms until either
$P=\emptyset$ or $|H|=p$. The procedure will stop in at most
$|J_{10}|+1$ steps for $\lambda^{i,J}_{r-1}$ and
$|J_{10}|$ steps for $\gamma^{i,J}_{r-1}$ because $|P| \leq |J_{10}|+1$ for
$\lambda^{i,J}_{r-1}$ and $|P| \leq |J_{10}|$ for $\gamma^{i,J}_{r-1}$
and on every step $|P|$ decreases at least by one.
Also, inequality $|T| \leq j$ is satisfied for the kernels of terms of
representations (\arabic{lambdaRepresentation}) and
(\arabic{gammaRepresentation}) and is preserved under transformation
from lemma~\ref{Transformation}, therefore the kernels of terms in the
right hand side of representation (\arabic{BRepresentation1})
with $P=\emptyset$ will be
${\cal K}^{M}_{a,b}$ or ${\cal B}^{\emptyset,M}_{a,0}$.\\
\indent
Conditions (\arabic{1StepIndices}) will
be satisfied because according to lemma~\ref{Transformation} after every
application of $Y_{i,\zeta}(z)$ to ${\widetilde g}(\zeta)$ numbers $b$ and
$|P|$ decrease by one with all other multiindices unchanged, therefore
decreasing $h\left({\cal B}^{P(H),M}_{a,b}\right)$ by one.\\
\indent
For each term of representation (\arabic{BRepresentation1}) with
$|H| = i$ we repeatedly apply lemma~\ref{Smoothness} or
lemma~\ref{LogSmoothness} $p-i$ times to obtain the following representation
$$D_1 \circ \dots \circ D_{p-i}
\left( \int_{U(\epsilon) \times [0,1]}
\left[ \{ Y_{\zeta}(z) \}^{H} {\widetilde g}(\zeta) \right]
{\cal B}^{P(H),M}_{a, b}(\zeta, z) dt \right)
\eqno(\arabic{equation})
\newcounter{BRepresentation2}
\setcounter{BRepresentation2}{\value{equation}}
\addtocounter{equation}{1}$$
$$= \sum_{r=i}^p \sum_{T, I} \int_{U(\epsilon) \times [0,1]}
c_{\left\{ T, I \right\} }(\zeta, z, t)
\left[ \{ Y_{\zeta}(z), D_{\zeta} \}^r {\widetilde g}(\zeta) \right]
{\cal B}^{T, I}_{d,j}(\zeta, z) dt,$$
\noindent
where $\{ Y_{\zeta}(z), D_{\zeta} \}^r$ denotes
a composition of differentiations $Y_{j, \zeta}(z)$ and
$D_{\zeta}$ applied $r$ times and indices $d,j$ and multiindices
$T, I$ of representation (\arabic{BRepresentation2}) satisfy
$$\begin{array}{llll}
k\left({\cal B}^{T, I}_{d,j}\right)+h\left({\cal B}^{T, I}_{d,j}\right)
-l\left({\cal B}^{T, I}_{d,j}\right)-s\left({\cal B}^{T, I}_{d,j}\right)
\vspace{0.1in}\\
\leq k\left({\cal B}^{P,M}_{a,b}\right)+h\left({\cal B}^{P,M}_{a,b}\right)
-l\left({\cal B}^{P,M}_{a,b}\right)-s\left({\cal B}^{P,M}_{a,b}\right)-i,
\vspace{0.1in}\\
k\left({\cal B}^{T, I}_{d,j}\right)+2h\left({\cal B}^{T, I}_{d,j}\right)
-2l\left({\cal B}^{T, I}_{d,j}\right)-2s\left({\cal B}^{T, I}_{d,j}\right)
\vspace{0.1in}\\
\leq k\left({\cal B}^{P,M}_{a,b}\right)+2h\left({\cal B}^{P,M}_{a,b}\right)
-2l\left({\cal B}^{P,M}_{a,b}\right)-2s\left({\cal B}^{P,M}_{a,b}\right)-2i.
\end{array}
\eqno(\arabic{equation})
\newcounter{2StepIndices}
\setcounter{2StepIndices}{\value{equation}}
\addtocounter{equation}{1}$$
\indent
Finally, we apply operator $D_1 \circ \dots \circ D_i$ to
each term of the representation (\arabic{BRepresentation2}) with
$r \leq p$ by differentiating the kernel and obtain
$$D_1 \circ \dots \circ D_i
\left( \int_{U(\epsilon) \times [0,1]}
\left[ \{ Y_{\zeta}(z), D_{\zeta} \}^r {\widetilde g}(\zeta) \right]
{\cal B}^{T, I}_{d,j}(\zeta, z) dt \right)
\eqno(\arabic{equation})
\newcounter{BRepresentation3}
\setcounter{BRepresentation3}{\value{equation}}
\addtocounter{equation}{1}$$
$$= \sum_{r \leq p} \sum_{P, M} \int_{U(\epsilon) \times [0,1]}
c_{\left\{ P, M \right\} }(\zeta, z, t)
\left[ \{ Y_{\zeta}(z), D_{\zeta} \}^r {\widetilde g}(\zeta) \right]
{\cal B}^{P, M}_{a, b}(\zeta, z) dt,$$
\noindent
with indices $a, b$ and multiindices $P, M$ satisfying
$$\begin{array}{llll}
k\left({\cal B}^{P,M}_{a,b}\right)+h\left({\cal B}^{P,M}_{a,b}\right)
-l\left({\cal B}^{P,M}_{a,b}\right)-s\left({\cal B}^{P,M}_{a,b}\right)
\vspace{0.1in}\\
\leq k\left({\cal B}^{T,I}_{d,j}\right)+h\left({\cal B}^{T,I}_{d,j}\right)
-l\left({\cal B}^{T,I}_{d,j}\right)-s\left({\cal B}^{T,I}_{d,j}\right),
\vspace{0.1in}\\
k\left({\cal B}^{P,M}_{a,b}\right)+2h\left({\cal B}^{P,M}_{a,b}\right)
-2l\left({\cal B}^{P,M}_{a,b}\right)-2s\left({\cal B}^{P,M}_{a,b}\right)
\vspace{0.1in}\\
\leq k\left({\cal B}^{T,I}_{d,j}\right)+2h\left({\cal B}^{T,I}_{d,j}\right)
-2l\left({\cal B}^{T,I}_{d,j}\right)-2s\left({\cal B}^{T,I}_{d,j}\right)
\end{array}
\eqno(\arabic{equation})
\newcounter{3StepIndices}
\setcounter{3StepIndices}{\value{equation}}
\addtocounter{equation}{1}$$
\indent
Conditions (\arabic{3StepIndices}) are satisfied because
after each application of $D_l$ to ${\cal B}^{T, I}_{d,j}(\zeta, z)$
either $k\left({\cal B}^{T,I}_{d,j}\right)$ increases by one
or $h\left({\cal B}^{T,I}_{d,j}\right)$ increases
by one or $h \left({\cal B}^{T,I}_{d,j} \right)$ increases by two and
$k\left({\cal B}^{T,I}_{d,j}\right)$ decreases by three.\\
\indent
From (\arabic{BRepresentation1}), (\arabic{BRepresentation2}) and
(\arabic{BRepresentation3}) we conclude that in order to prove the
statement of the proposition it suffices to prove that
$$\left\| \int_{U(\epsilon) \times [0,1]} c(\zeta, z, t)
\left[ \{ Y_{\zeta}(z), D_{\zeta} \}^p {\widetilde g}(\zeta) \right]
{\cal B}^{P, M}_{a, b}(\zeta, z) dt \right\|_{ \Gamma^{p,\alpha +1}(U) } \leq
C \cdot \| g \|_{ \Gamma^{p, \alpha}(U) }
\eqno(\arabic{equation})
\newcounter{FinalEstimate}
\setcounter{FinalEstimate}{\value{equation}}
\addtocounter{equation}{1}$$
for the kernel ${\cal B}^{P, M}_{a, b}$ obtained from $\lambda^{i,J}_{r-1}$
after described above application of lemmas \ref{Smoothness},
\ref{Transformation} and \ref{LogSmoothness}.\\
\indent
We will prove (\arabic{FinalEstimate}) as a corollary of the lemma below.

\begin{lemma}\label{GammaAlpha}
Let $0< \alpha <1$, $g \in \Gamma^{\alpha}
\left( \{ \rho \}, {\widetilde U}(\epsilon_0) \right)$
be a function with compact support and let indices $a,b$
and multiindices $P,M$ satisfy conditions
$$\begin{array}{ll}
k\left({\cal B}^{P,M}_{a,b}\right)+h\left({\cal B}^{P,M}_{a,b}\right)
-l\left({\cal B}^{P,M}_{a,b}\right)-s\left({\cal B}^{P,M}_{a,b}\right)
\leq 2n-m-2, \vspace{0.1in}\\
k\left({\cal B}^{P,M}_{a,b}\right)+2h\left({\cal B}^{P,M}_{a,b}\right)
-2l\left({\cal B}^{P,M}_{a,b}\right)-2s\left({\cal B}^{P,M}_{a,b}\right)
\leq 2n-1
\end{array}
\eqno(\arabic{equation})
\newcounter{FinalIndices}
\setcounter{FinalIndices}{\value{equation}}
\addtocounter{equation}{1}$$
for $b \geq 1$ and
$$k\left({\cal B}^{\emptyset,M}_{a,0}\right)
-l\left({\cal B}^{\emptyset,M}_{a,0}\right)
-s\left({\cal B}^{\emptyset,M}_{a,0}\right)
\leq 2n-m-3
\eqno(\arabic{equation})
\newcounter{LogFinalIndices}
\setcounter{LogFinalIndices}{\value{equation}}
\addtocounter{equation}{1}$$
for $|P|=0$ and $b=0$.\\
\indent
Then 
$$f_{\epsilon}(z):= \left( \int_{U(\epsilon) \times [0,1]} c(\zeta, z, t)
g(\zeta) {\cal B}^{P, M}_{a, b}(\zeta, z) dt \right)
\in \Gamma^{\alpha+1}(U)$$
\noindent
for $\epsilon < \epsilon_0$ and
$$\left\| f_{\epsilon} \right\|_{ \Gamma^{\alpha+1}(U) } \leq
C \cdot \| g \|_{ \Gamma^{\alpha}(U) }$$
\noindent
with $C$ independent of $g$ and $\epsilon$.
\end{lemma}

\indent
{\bf Proof.}\\
\indent
At first we will prove the inclusion $f_{\epsilon} \in
\Lambda^{\frac{\alpha+1}{2} }(U)$.
For a fixed point $w \in U$ and arbitrary $z \in U$ we denote $\delta = |z-w|$
and $c = 1+ \max_{i, j, z \in U}\left\{ |Q^{(k)}_j(z)|\right\}$ and introduce
the following neighborhood
$$V(w,z)= \left\{ \zeta \in U :
|\zeta - w|^2 + \sum_{i=1}^m |\rho_i(\zeta)|
+ \sum_{i=1}^m \left| \mbox{Im}F^{(i)}(\zeta,w)\right| \leq
4cn^2\delta \right\},$$
\noindent
containing $z$ and $w$ and such that the estimates
$$\left| \Phi(\zeta,w) \right|,\left| \Phi(\zeta,z) \right|
\geq C \cdot \delta^2$$
hold for $\zeta \notin V(w,z)$ with constant $C>0$ independent of $\delta$
for $\delta$ small enough.\\
\indent
Denoting $V(\epsilon) = V(w,z) \cap U(\epsilon)$ we represent
$f_{\epsilon}(z)$ as follows
$$f_{\epsilon}(z) = g(w) \cdot \int_{U(\epsilon) \times [0,1]} c(\zeta, z, t)
{\cal B}^{P, M}_{a, b}(\zeta, z) dt
\eqno(\arabic{equation})
\newcounter{fRepresentation1}
\setcounter{fRepresentation1}{\value{equation}}
\addtocounter{equation}{1}$$
$$+ \int_{V(\epsilon) \times [0,1]} \left(g(\zeta)-g(w)\right)
c(\zeta, z, t) {\cal B}^{P, M}_{a, b}(\zeta, z) dt$$
$$+ \int_{\left( U(\epsilon) \setminus V(\epsilon) \right) \times [0,1]}
\left(g(\zeta)-g(w)\right) c(\zeta, z, t)
{\cal B}^{P, M}_{a, b}(\zeta, z) dt.$$
\indent
Applying lemmas \ref{Smoothness}, \ref{Transformation}
and \ref{LogSmoothness} as earlier in the proof of
proposition~\ref{REstimate} to the first
term of the right hand side of (\arabic{fRepresentation1}) we obtain
$$D \left( g(w) \int_{U(\epsilon) \times [0,1]} c(\zeta, z, t) \cdot
{\cal B}^{P,M}_{a,b}(\zeta, z) dt \right)$$
$$=g(w) \cdot \sum_{T,I,d,j} \int_{U(\epsilon) \times [0,1]}
c_{ \{ T,I,d,j \}}(\zeta, z, t)
\cdot {\cal B}^{T,I}_{d,j}(\zeta, z) dt$$
\noindent
with $c_{ \{ T,I,d,j \}}(\zeta, z, t)
=c_{ \{ T,I,d,j \}}(\zeta,z,\theta(\zeta),t)$ and
$$c_{ \{ T,I,d,j \}}(\zeta,z,\theta,t)
\in C^{\infty} \left( {\widetilde U}_{\zeta} \times {\widetilde U}_{z}
\times \S^{n-1} \times [0,1] \right)$$
\noindent
and indices $d,j$ and multiindices $T, I_1, I_2, I_3$ satisfying
(\arabic{FinalIndices}). Then applying formula
$$\left. d\rho_i \right|_{U(\epsilon)} = \epsilon d\theta_i$$
\noindent
and lemma~\ref{Integral} we obtain
$$\left| g(w) \cdot D \left( \int_{U(\epsilon) \times [0,1]} c(\zeta, z, t)
{\cal B}^{P, M}_{a, b}(\zeta, z) dt \right) \right|$$
$$= \| g \|_{ \Lambda^{0}(U) } \cdot {\cal O} \left(
\int_{U(\epsilon) \times [0,1]}
\frac{ \overbrace{\wedge d_{\zeta}\mbox{Re}F}^{s({\cal B})}
\overbrace{\wedge_{i} d\theta_i(\zeta)}^{m-1-s({\cal B})}
\wedge d\sigma_{2n-m}(\zeta)} {|\zeta - z|^{k({\cal B})} \cdot
|\Phi(\zeta,z)|^{\frac{h({\cal B})-l({\cal B})}{2}} } \right)$$
$$= \| g \|_{ \Lambda^{0}(U) } \cdot
{\cal O} \left({\cal I}_1\left\{ 0,k({\cal B}), h({\cal B})-l({\cal B}),
s({\cal B}) \right\}
\left( \epsilon, 1 \right) \right)
= \| g \|_{ \Lambda^{0}(U) } \cdot {\cal O}(1),$$
\noindent
which shows that the first term of the right hand side of
(\arabic{fRepresentation1}) is in $\Lambda^1(U)$.\\
\indent
For the second term of the right hand side of (\arabic{fRepresentation1})
we have
$$\left| \int_{V(\epsilon) \times [0,1]} \left(g(\zeta)-g(w)\right)
c(\zeta, z, t) {\cal B}^{P, M}_{a, b}(\zeta, z) dt \right|$$
$$= \| g \|_{ \Gamma^{\alpha}(U) } \cdot
\delta^{ \frac{a}{2} } \cdot {\cal O} \left(
\int_{V(\epsilon) \times [0,1]}
\frac{ \overbrace{\wedge d_{\zeta}\mbox{Re}F}^{s({\cal B})}
\overbrace{\wedge_{i} d\theta_i(\zeta)}^{m-1-s({\cal B})}
\wedge d\sigma_{2n-m}(\zeta)} {|\zeta - z|^{k({\cal B})} \cdot
|\Phi(\zeta,z)|^{\frac{h({\cal B})-l({\cal B})}{2}} } \right)$$
$$= \| g \|_{ \Gamma^{\alpha}(U) } \cdot \delta^{ \frac{\alpha}{2} } \cdot
{\cal O} \left({\cal I}_1\left\{ 0,k({\cal B}), h({\cal B})-l({\cal B}),
s({\cal B})\right\}
\left( \epsilon, \sqrt{\delta} \right) \right)
= \| g \|_{ \Gamma^{\alpha}(U) } \cdot
{\cal O} \left( \delta^{ \frac{\alpha+1}{2} } \right),$$
\noindent
where we used lemma~\ref{Integral} and the estimate
$$|g(\zeta) - g(w)| = \| g \|_{ \Gamma^{\alpha}(U) }
\cdot {\cal O}\left( \delta^{ \frac{\alpha}{2} } \right)$$
\noindent
for $\zeta \in V(w,z)$.\\
\indent
For the third term of the right hand side of (\arabic{fRepresentation1})
we have
$$\left| \int_{\left( U(\epsilon) \setminus V(\epsilon) \right) \times [0,1]}
\left(g(\zeta)-g(w)\right)
c(\zeta, z, t) {\cal B}^{P, M}_{a, b}(\zeta, z) dt \right.$$
$$\left.- \int_{\left( U(\epsilon) \setminus V(\epsilon) \right) \times [0,1]}
\left(g(\zeta)-g(w)\right)
c(\zeta, w, t){\cal B}^{P, M}_{a, b}(\zeta, w) dt \right|$$
$$\leq \left|
\int_{\left( U(\epsilon) \setminus V(\epsilon) \right) \times [0,1]}
\left(g(\zeta)-g(w)\right)
c(\zeta, z, t) \left[ {\cal B}^{P, M}_{a, b}(\zeta, z)
- {\cal B}^{P, M}_{a, b}(\zeta, w) \right] dt\right|$$
$$+ \left| \int_{\left( U(\epsilon) \setminus V(\epsilon) \right) \times [0,1]}
\left(g(\zeta)-g(w)\right)
\left[ c(\zeta, z, t) - c(\zeta.w) \right]
{\cal B}^{P, M}_{a, b}(\zeta, w) dt\right|$$
$$= \| g \|_{ \Gamma^{\alpha}(U) } \cdot \delta \cdot {\cal O} \left[
{\cal I}_2\left\{ \alpha, k({\cal B})+1, h({\cal B})-l({\cal B}),
s({\cal B}) \right\}
\left( \epsilon, \sqrt{\delta} \right) \right.$$
$$+ {\cal I}_2\left\{ \alpha, k({\cal B}),
h({\cal B})-l({\cal B})+1, s({\cal B}) \right\}
\left( \epsilon,\sqrt{\delta} \right)$$
$$\left.+ {\cal I}_2\left\{ \alpha, k({\cal B}),
h({\cal B})-l({\cal B}), s({\cal B}) \right\}
\left( \epsilon,\sqrt{\delta} \right) \right] =
\| g \|_{ \Gamma^{\alpha}(U) } \cdot
{\cal O}\left( \delta^{ \frac{\alpha+1}{2} } \right),$$
\noindent
where we again used lemma~\ref{Integral}.\\
\indent
Representation (\arabic{fRepresentation1}) together with the estimates
above show that 
$$\left\| f_{\epsilon} \right\|_{ \Lambda^{ \frac{\alpha+1}{2} }(U) } \leq
C \cdot \| g \|_{ \Gamma^{\alpha}(U) }$$
uniformly with respect to $\epsilon$.\\
\indent
To complete the proof of the lemma we have to prove that
$$\left\| D^c f_{\epsilon} \right\|_{ \Lambda^{\alpha}_c(U) } \leq
C \cdot \| g \|_{ \Gamma^{\alpha}(U) }$$
\noindent
where differentiation $D^c \in T^c({{\bold M}})_{z}$ and
$$\left\| h \right\|_{ \Lambda^{\alpha}_c(U) }
= \sup \left\{ \left| h(x(\cdot)) \right|_{\Lambda^{\alpha}([0,1])} \right\}$$
\noindent
with the {\it sup} taken over all curves
$x : [0,1] \rightarrow {\bold M}$ such that
$$\begin{array}{ll}
(i)\hspace{0.1in} |x^{\prime}(t)|, |x^{\prime\prime}(t)| \leq 1,
\vspace{0.1in}\\
(ii)\hspace{0.1in}x^{\prime}(t) \in T^{c}({\bold M}).
\end{array}$$
\indent
To prove this estimate we use the following representation
$$D^c f_{\epsilon}(z)
= g(z) \cdot  D^c \left( \int_{U(\epsilon) \times [0,1]}
c(\zeta, z, t) {\cal B}^{P, M}_{a, b}(\zeta, z) dt \right)
\eqno(\arabic{equation})
\newcounter{DcfRepresentation1}
\setcounter{DcfRepresentation1}{\value{equation}}
\addtocounter{equation}{1}$$
$$+ \left( \int_{U(\epsilon) \times [0,1]} \left(g(\zeta)-g(z)\right)
D^c \left[c(\zeta, z, t) {\cal B}^{P, M}_{a, b}(\zeta, z)\right] dt \right).$$
\indent
Applying as for the first term of (\arabic{fRepresentation1}) lemmas
\ref{Transformation}, \ref{Smoothness} and \ref{LogSmoothness} and using
lemma~\ref{Integral} we conclude that
$$D^c \left(\int_{U(\epsilon) \times [0,1]} c(\zeta, z, t)
{\cal B}^{P, M}_{a, b}(\zeta, z) dt \right) \in \Lambda^1(U)$$
and therefore the first term of the right hand side of
(\arabic{DcfRepresentation1}) is in $\Lambda^{\alpha}(U)$.\\
\indent
To estimate the second term of the right hand side of (\arabic{DcfRepresentation1})
we consider $z, w \in U$ and introduce the neighborhoods
$$W(w,z)= \left\{ \zeta \in U :
|\zeta - w|^2 + \sum_{i=1}^m |\rho_i(\zeta)|
+ \sum_{i=1}^m \left| \mbox{Im}F^{(i)}(\zeta,w)\right| \leq
4cn^2\delta^2 \right\},$$
$$W^{\prime}(w,z)= \left\{ \zeta \in U :
|\zeta - w|^2 + \sum_{i=1}^m |\rho_i(\zeta)|
+ \sum_{i=1}^m \left| \mbox{Im}F^{(i)}(\zeta,w)\right| \leq
16cn^2\delta^2 \right\},$$
$$W(\epsilon) = W(w,z) \cap U(\epsilon),$$
$$W^{\prime}(\epsilon) = W^{\prime}(w,z) \cap U(\epsilon)$$
and a function $\phi(\zeta) \in C^{\infty}({\widetilde U})$ such that
$0 \leq \phi(\zeta) \leq 1$, $\phi \equiv 1$ on $W(w,z)$,
$\phi \equiv 0$ on ${\widetilde U} \setminus W^{\prime}(w,z)$,
$\left| \mbox{grad}\phi(\zeta) \right| \leq 2/\delta^2$ and
$\left| \mbox{grad}_c\phi(\zeta) \right| \leq 2/\delta$.\\
\indent
Then we consider the following representation of the second term of the
right hand side of (\arabic{DcfRepresentation1})
$$\int_{U(\epsilon) \times [0,1]} \left(g(\zeta)-g(z)\right)
D^c \left[c(\zeta, z, t) {\cal B}^{P, M}_{a, b}(\zeta, z)\right] dt
\eqno(\arabic{equation})
\newcounter{DcfIntegral3}
\setcounter{DcfIntegral3}{\value{equation}}
\addtocounter{equation}{1}$$
$$-\int_{U(\epsilon) \times [0,1]} \left(g(\zeta)-g(w)\right)
D^c \left[ c(\zeta, w, t) {\cal B}^{P, M}_{a, b}(\zeta, w) \right]dt$$
$$=\int_{U(\epsilon) \times [0,1]} \left(g(\zeta)-g(z)\right) \phi(\zeta)
D^c \left[ c(\zeta, z, t) {\cal B}^{P, M}_{a, b}(\zeta, z)\right] dt$$
$$-\int_{U(\epsilon) \times [0,1]} \left(g(\zeta)-g(w)\right) \phi(\zeta)
D^c \left[ c(\zeta, w, t) {\cal B}^{P, M}_{a, b}(\zeta, w)\right] dt$$
$$+\int_{ U(\epsilon) \times [0,1] }
\left[ \left(g(\zeta)-g(z)\right) \left( 1-\phi(\zeta) \right)
D^c \left[ c(\zeta, z, t) {\cal B}^{P, M}_{a, b}(\zeta, z)\right] \right.$$
$$\left. - \left(g(\zeta)-g(w)\right) \left( 1-\phi(\zeta) \right)
D^c \left[c(\zeta, w, t) {\cal B}^{P, M}_{a, b}(\zeta, w)\right] \right] dt.$$
\indent
To estimate integrals on the right hand side of (\arabic{DcfIntegral3})
we use the estimate
$$D^c \left[ F^{(k)} (\zeta,z) \right] = {\cal O}\left(|\zeta -z|\right)
\eqno(\arabic{equation})
\newcounter{DcF}
\setcounter{DcF}{\value{equation}}
\addtocounter{equation}{1}$$
for $k=1, \dots, m$ and obtain a representation
$$D^c \left[ {c(\zeta, z, t) \cal B}^{P, M}_{a, b}(\zeta, z)\right] =
\sum_{T,I,d,j} c_{ \{ T, I, d,j \} }(\zeta, z, t)
\cdot {\cal B}^{T, I}_{d,j}(\zeta, z)
\eqno(\arabic{equation})
\newcounter{DcB}
\setcounter{DcB}{\value{equation}}
\addtocounter{equation}{1}$$
with $c_{ \{ T,I,d,j \}}(\zeta, z, t)
=c_{ \{ T,I,d,j \}}(\zeta,z,\theta(\zeta),t)$ and
$$c_{ \{ T,I,d,j \}}(\zeta,z,\theta,t)
\in C^{\infty} \left( {\widetilde U}_{\zeta} \times {\widetilde U}_{z}
\times \S^{n-1} \times [0,1] \right),$$
\noindent
and indices and multiindices in
the right hand side of (\arabic{DcB}) satisfying
$$\begin{array}{ll}
k\left({\cal B}^{T,I}_{d,j}\right)+h\left({\cal B}^{T,I}_{d,j}\right)
-l\left({\cal B}^{T,I}_{d,j}\right)-s\left({\cal B}^{T,I}_{d,j}\right)
\leq 2n-m-1, \vspace{0.1in}\\
k\left({\cal B}^{T,I}_{d,j}\right)+2h\left({\cal B}^{T,I}_{d,j}\right)
-2l\left({\cal B}^{T,I}_{d,j}\right)-2s\left({\cal B}^{T,I}_{d,j}\right)
\leq 2n.
\end{array}
\eqno(\arabic{equation})
\newcounter{DcBIndices}
\setcounter{DcBIndices}{\value{equation}}
\addtocounter{equation}{1}$$
\indent
Then for the first two integrals of the right hand side of
(\arabic{DcfIntegral3}) using representation (\arabic{DcB}) with conditions
(\arabic{DcBIndices}) and lemma~\ref{Integral} we obtain
$$\left| \int_{U(\epsilon) \times [0,1]} \left(g(\zeta)-g(z)\right)
\phi(\zeta) c(\zeta, z, t) {\cal B}^{T,I}_{d,j}(\zeta, z) dt \right|$$
$$={\cal O} \left( \left|
\int_{W^{\prime}(\epsilon) \times [0,1]} \left(g(\zeta)-g(z)\right)
c(\zeta, z, t) {\cal B}^{T,I}_{d,j}(\zeta, z) dt \right| \right)$$
$$= \| g \|_{ \Gamma^{\alpha}(U) } \cdot {\cal O} \left(
\int_{W(\epsilon)}\frac{ \overbrace{\wedge d_{\zeta}\mbox{Re}F}^{s({\cal B})}
\overbrace{\wedge_{i} d\theta_i(\zeta)}^{m-1-s({\cal B})}
\wedge d\sigma_{2n-m}(\zeta)} {|\zeta - z|^{k({\cal B})} \cdot
|\Phi(\zeta,z)|^{\frac{h({\cal B})-l({\cal B})}{2}-\frac{\alpha}{4}} }
\right)$$
$$= \| g \|_{ \Gamma^{\alpha}(U) } \cdot {\cal O}
\left({\cal I}_1\left\{ \alpha, k({\cal B}), h({\cal B})-l({\cal B}),
s({\cal B}) \right\}
\left( \epsilon, \delta \right) \right)
= \| g \|_{ \Gamma^{\alpha}(U) } \cdot
{\cal O} \left( \delta^{ \alpha } \right).$$
\indent
To estimate the third integral of the right hand side of
(\arabic{DcfIntegral3}) we represent it as
$$\int_{ U(\epsilon) \times [0,1]}
\left(g(\zeta)-g(z)\right) \left( 1-\phi(\zeta) \right) 
D^c \left[c(\zeta, z, t){\cal B}^{P, M}_{a, b}(\zeta, z)\right] dt
\eqno(\arabic{equation})
\newcounter{DcfIntegral4}
\setcounter{DcfIntegral4}{\value{equation}}
\addtocounter{equation}{1}$$
$$- \int_{ U(\epsilon) \times [0,1]}
\left(g(\zeta)-g(w)\right) \left( 1-\phi(\zeta) \right)
D^c \left[c(\zeta, w, t){\cal B}^{P, M}_{a, b}(\zeta, w)\right] dt$$
$$= \int_{ U(\epsilon) \times [0,1]}
\left(g(\zeta)-g(z)\right) \left( 1-\phi(\zeta) \right)$$
$$\times\left( D^c \left[c(\zeta, z, t){\cal B}^{P, M}_{a, b}(\zeta, z)\right]
- D^c \left[c(\zeta, w, t){\cal B}^{P, M}_{a, b}(\zeta, w) \right] \right)dt$$
$$+\left[g(w)-g(z)\right] \cdot \int_{ U(\epsilon) \times [0,1]}
\left( 1-\phi(\zeta) \right)
D^c \left[c(\zeta, w, t){\cal B}^{P, M}_{a, b}(\zeta, w)\right] dt.$$
\indent
To estimate the first integral of the right hand side of
(\arabic{DcfIntegral4}) we use representation (\arabic{DcB}) with
kernels satisfying (\arabic{DcBIndices}) and
estimate
$$\left| F^{(k)}(\zeta, z) - F^{(k)}(\zeta, w) \right|
= {\cal O} \left( \delta \cdot |\zeta - z| + \delta^2 \right)
\eqno(\arabic{equation})
\newcounter{Fz-Fw}
\setcounter{Fz-Fw}{\value{equation}}
\addtocounter{equation}{1}$$
for $k = 1, \dots, m$ and $\zeta \in U \setminus W(\epsilon)$.\\
\indent
Then we obtain
$$\left| \int_{ U(\epsilon) \times [0,1]}
\left(g(\zeta)-g(z)\right) \left( 1-\phi(\zeta) \right) \right.$$
$$\left.\times
\left( D^c \left[c(\zeta, z, t){\cal B}^{P, M}_{a, b}(\zeta, z)\right]
- D^c \left[c(\zeta, w, t){\cal B}^{P, M}_{a, b}(\zeta, w) \right] \right)dt
\right|$$
$$=\| g \|_{ \Gamma^{\alpha}(U) } \cdot \left[\delta \cdot {\cal O}
\left({\cal I}_2\left\{ \alpha, k({\cal B}), h({\cal B})-l({\cal B}),
s({\cal B})\right\} \left( \epsilon,\delta \right) \right) \right.$$
$$+ \delta \cdot {\cal O}
\left({\cal I}_2\left\{ \alpha, k({\cal B})+1, h({\cal B})-l({\cal B}),
s({\cal B})\right\}
\left( \epsilon,\delta \right) \right)$$
$$+ \delta \cdot {\cal O}
\left({\cal I}_2\left\{ \alpha, k({\cal B})-1, h({\cal B})-l({\cal B})+1,
s({\cal B})\right\} \left( \epsilon,\delta \right) \right)$$
$$\left.+ \delta^2 \cdot {\cal O}
\left({\cal I}_2\left\{ \alpha, k({\cal B}), h({\cal B})-l({\cal B})+1,
s({\cal B})\right\} \left( \epsilon,\delta \right) \right) \right]
= \| g \|_{ \Gamma^{\alpha}(U) } \cdot
{\cal O}\left( \delta^{\alpha} \right).$$
\indent
To obtain necessary estimate for the second integral of the
right hand side of (\arabic{DcfIntegral4}) it suffices to prove the following
estimate
$$\left| \int_{ U(\epsilon) \times [0,1]}
\left( 1-\phi(\zeta) \right)
D^c \left[c(\zeta, w, t){\cal B}^{P, M}_{a, b}(\zeta, w)\right]dt \right|
= {\cal O}(1).
\eqno(\arabic{equation})
\newcounter{DcfEstimate}
\setcounter{DcfEstimate}{\value{equation}}
\addtocounter{equation}{1}$$
\indent
To obtain this estimate for kernels ${\cal B}^{P, M}_{a, b}$ with
$|P| = 0$ and $b \geq 1$ (i.e. for ${\cal K}^{M}_{a, b}$) we use
representation from lemma~\ref{Smoothness}
$$D^c \int_{ U(\epsilon) \times [0,1]} \left( 1-\phi(\zeta) \right)
c(\zeta, w, t){\cal K}^{M}_{a, b}(\zeta, w)dt
\eqno(\arabic{equation})
\newcounter{DcRepresentation1}
\setcounter{DcRepresentation1}{\value{equation}}
\addtocounter{equation}{1}$$
$$= \int_{U(\epsilon) \times [0,1]} \left( 1-\phi(\zeta) \right)
\left[ D^c c(\zeta, w, t) \right] \cdot {\cal K}^{M}_{a,b}(\zeta, w)dt$$
$$+ \sum_{S,d,j} \int_{U(\epsilon) \times [0,1]}
c_{ \{ S,d,j \}}(\zeta, w, t) \cdot \left( 1-\phi(\zeta) \right)
\cdot {\cal K}^{S}_{d,j}(\zeta, w) dt$$
$$+ \int_{U(\epsilon) \times [0,1]}
D^c_{\zeta} \left( 1-\phi(\zeta) \right) \cdot c(\zeta, w, t)
{\cal K}^{M}_{a,b}(\zeta, w) dt$$
$$+ \sum_{k=1}^{m} \sum_{L} \int_{U(\epsilon) \times [0,1]}
\left[ Y_{k,\zeta}(w) \left( 1-\phi(\zeta) \right) \right]
\cdot c_{ \{ k, L \}}(\zeta, w, t) \cdot {\cal K}^{L}_{a,b}(\zeta, w) dt$$
with indices of the kernels of all the terms except the last term satisfying
(\arabic{FinalIndices}). For the kernels of the last term of the
right hand side of (\arabic{DcRepresentation1}) according to 
(\arabic{ImIntegral1}) and (\arabic{ImIntegral2}) we have
$$k\left({\cal K}^{L}_{a,b}\right) = k\left({\cal K}^{M}_{a,b}\right)-1.
\eqno(\arabic{equation})
\newcounter{DcIndices}
\setcounter{DcIndices}{\value{equation}}
\addtocounter{equation}{1}$$
\indent
For the first two terms of the right hand side of (\arabic{DcRepresentation1})
we obtain estimate (\arabic{DcfEstimate}) using lemma~\ref{Integral} and
conditions (\arabic{FinalIndices}). For the third term of the right hand side
of (\arabic{DcRepresentation1}) we obtain the following estimate
$$\left| \int_{U(\epsilon) \times [0,1]}
D^c_{\zeta} \left( 1-\phi(\zeta) \right) \cdot c(\zeta, w, t)
{\cal K}^{M}_{a,b}(\zeta, w) dt \right|$$
$$= {\delta}^{-1} \cdot {\cal O} \left(
{\cal I}_1\left\{ 0, k({\cal K}), h({\cal K})-l({\cal K}), s({\cal K}) \right\}
\left( \epsilon, \delta \right) \right)
= {\cal O}(1),$$
\noindent
where we used properties of the function $\phi$.\\
\indent
We obtain the same estimate for the integrals of the fourth term of the right
hand side of (\arabic{DcRepresentation1}) if we use the inequality
$$\left| \zeta - w \right| = {\cal O}( \delta )\hspace{0.1in}\mbox{for}
\hspace{0.1in}\zeta \in W^{\prime}(w,z),$$
\noindent
property (\arabic{DcIndices}) and the estimate
$\left| \mbox{grad}\phi(\zeta) \right| \leq 2/\delta^2$.\\
\indent
Proof of the same estimates for kernels ${\cal B}^{\emptyset,M}_{a,0}$ is
analogous to the proof for ${\cal K}^{M}_{a,b}$ but uses condition
(\arabic{LogFinalIndices}) instead of conditions (\arabic{FinalIndices}).\\
\indent
For a general kernel ${\cal B}^{P, M}_{a, b}$ with $|P| \neq 0$ we use
representation from lemma~\ref{Transformation} and obtain
$$D^c \int_{U(\epsilon) \times [0,1]}
\left( 1-\phi(\zeta) \right) c(\zeta, w, t)
\cdot {\cal B}^{P,M}_{a,b}(\zeta, w) dt
\eqno(\arabic{equation})
\newcounter{DcRepresentation2}
\setcounter{DcRepresentation2}{\value{equation}}
\addtocounter{equation}{1}$$
$$= \int_{U(\epsilon) \times [0,1]}
\left[ Y_{k,\zeta}(z) \left( 1-\phi(\zeta) \right) \right] \cdot
D^c \left[ c(\zeta, w, t) {\cal B}^{T,I}_{a,b-1}(\zeta, w) \right]dt$$
$$+ \sum_{ \left\{ L, S, d, j \right\} }
D^c \int_{U(\epsilon) \times [0,1]} \left( 1-\phi(\zeta) \right)
c_{\left\{ L, S, d, j \right\} }(\zeta, w, t)
\cdot {\cal B}^{L,S}_{d,j}(\zeta, w)dt,$$
with $|T| = |P|-1$, indices of the kernels ${\cal B}^{L,S}_{d,j}$ satisfying
(\arabic{FinalIndices}) and $|L| < |P|$. Therefore the problem
of estimating integral in (\arabic{DcfEstimate}) will be reduced to estimating
the same integral for a kernel with smaller $|P|$ if we prove estimate
(\arabic{DcfEstimate}) for the second term of the right hand side of
(\arabic{DcRepresentation2}). We obtain this estimate analogously to
the estimate of the fourth term of the right hand side of
(\arabic{DcRepresentation1}) with the use of estimate (\arabic{DcF}).\qed

\indent
In order to prove applicability of lemma~\ref{GammaAlpha} to the kernels
obtained from $\lambda^{i,J}_{r-1}$ and $\gamma^{i,J}_{r-1}$ after
applications of lemmas \ref{Smoothness}, \ref{Transformation}
and \ref{LogSmoothness} we have
to prove relations (\arabic{FinalIndices}) and (\arabic{LogFinalIndices})
for these kernels.
But according to lemmas \ref{Smoothness}, \ref{Transformation}
and \ref{LogSmoothness} expressions in the left hand sides of these
relations don't increase under transformations from these lemmas.
Also, as follows from lemma~\ref{Transformation}, relation
(\arabic{LogFinalIndices}) for kernels ${\cal B}^{\emptyset,M}_{a,0}$
is a corollary of the first relation from (\arabic{FinalIndices}) for
kernels ${\cal B}^{P,M}_{a,1}$ with $|P| = 1$.
Therefore it suffices to prove conditions (\arabic{FinalIndices}) for
the original kernels ${\cal B}^{T,I}_{d,j}(\zeta, z)$ satisfying conditions
(\arabic{KlambdaIndices}) and (\arabic{KgammaIndices}).\\
\indent
We notice that since $s\left({\cal B}^{T,I}_{d,j}\right) = 0$ for the
original kernels we will omit it in our calculations.\\
\indent
Second condition from (\arabic{FinalIndices}) is always satisfied for
the indices satisfying (\arabic{KlambdaIndices})
as can be seen from the inequality
$$l({\cal B})+\frac{1}{2} \left( 2n-k({\cal B})-2h({\cal B}) \right)
\eqno(\arabic{equation})
\newcounter{SecondIndices}
\setcounter{SecondIndices}{\value{equation}}
\addtocounter{equation}{1}$$
$$\geq \frac{1}{2} \left( |J_8| + m - |J_2| - |J_3| - |J_{11}|\right)
\geq \frac{1}{2},$$
where we used relations
$$|J_1| + |J_3| + |J_6|+r+m-n \geq 0,$$
$$\sum_{i=1}^6 |J_{i}| = n-r-1,$$
$$|J_2| + |J_3| + |J_{11}| \leq m-1$$
for the multiindices of $\lambda^{i,J}_{r-1}$.\\
\indent
Condition (\arabic{SecondIndices}) is also satisfied for
the indices satisfying (\arabic{KgammaIndices}) which follows
from the inequality
$$l({\cal B})+\frac{1}{2} \left( 2n-k({\cal B})-2h({\cal B}) \right)$$
$$\geq \frac{1}{2} \left( 1+|J_8| + m - |J_2| - |J_3| - |J_{11}|\right)
\geq \frac{1}{2},$$
where we used relations
$$|J_1| + |J_3| + |J_6|+r+m-n \geq 0,$$
$$\sum_{i=1}^6 |J_{i}| = n-r-1,$$
$$|J_2| + |J_3| + |J_{11}| \leq m$$
for the multiindices of $\gamma^{i,J}_{r-1}$.\\
\indent
First condition from (\arabic{FinalIndices}) is not satisfied
for all kernels ${\cal B}^{T,I}_{d,j}(\zeta, z)$. But in the lemma
below we prove that if it is not satisfied then the corresponding terms
of the integral formula for $R_r(\epsilon)$ do not survive under the
limit when $\epsilon \rightarrow 0$. This lemma is a simplified version
of the lemma 4 in \cite{P}.\\

\begin{lemma}\label{RDoomed}
If
$$k({\cal B})+h({\cal B})-l({\cal B}) \geq 2n-m-1
\hspace{0.1in}\mbox{and}\hspace{0.1in}s({\cal B})=0$$
\noindent
then
$$\left\| \int_{U(\epsilon) \times [0,1]}
{\widetilde g}(\zeta) c(\zeta,z,t) {\cal B}^{T,I}_{d,j}(\zeta, z) dt
\right\|_{L^{\infty}({\bold M})} =
{\cal O}(\sqrt{\epsilon} \cdot \log{\epsilon})
\cdot \| g \|_{L^{\infty} ({\bold M})}.$$
\end{lemma}

\indent
{\bf Proof.}\\
\indent
We use the inequality
$$2n-m+ l({\cal B}) - k({\cal B}) - h({\cal B})
\geq n-|J_1|-|J_2|-|J_7|-1 \geq 1,
\eqno(\arabic{equation})
\newcounter{FirstIndices}
\setcounter{FirstIndices}{\value{equation}}
\addtocounter{equation}{1}$$
which is a corollary of the definitions of $k({\cal B})$, $h({\cal B})$
and $l({\cal B})$ and equality
$$\sum_{i=1}^6 |J_{i}| = n-r-1.$$
\indent
From the condition of the lemma and the inequality (\arabic{FirstIndices})
we obtain
$$k({\cal B})+h({\cal B})-l({\cal B}) = 2n-m-1$$
\noindent
and
$$n-|J_1|-|J_2|-|J_7|-1 = 1,$$
which leads to
$$|J_1|+|J_2|=n-r-1,\hspace{0.1in}|J_3|=0,\hspace{0.1in}|J_6|=0,
\hspace{0.1in}|J_7|=r-1,$$
and hence to
$$l({\cal B}) \geq |J_1|+|J_2|+|J_3|+|J_6|+r+m-n = m-1 > 1.$$
\indent
Using lemma~\ref{Integral} in the estimate of the integral in lemma
we obtain
$$\left| \int_{U(\epsilon) \times [0,1]}
{\widetilde g}(\zeta) c(\zeta,z,t)
{\cal B}^{T,I}_{d,j}(\zeta, z) dt \right|$$
$$= \| g \|_{L^{\infty} ({\bold M})} \cdot
\epsilon^{l({\cal B})} \cdot {\cal O}
\left( {\cal I}_1\left\{ 0, k({\cal B}), h({\cal B}), 0\right\}
\left( \epsilon, 1 \right) \right)$$
$$= \| g \|_{L^{\infty} ({\bold M})} \cdot
\left\{ \begin{array}{ll}
\epsilon^{l({\cal B})} \cdot
{\cal O}\left(\epsilon^{2n-m-k({\cal B})-h({\cal B})}
\cdot (\log{\epsilon})^2 \right) & \mbox{if} \hspace{0.05in}
k({\cal B}) \geq 2n-2m,
\vspace{0.1in}\\
\epsilon^{l({\cal B})} \cdot
{\cal O}\left( \epsilon^{(2n-k({\cal B})-2h({\cal B}))/2}
\cdot \log{\epsilon}\right) & \mbox{if} \hspace{0.05in}
k({\cal B}) \leq 2n-2m-1.
\end{array} \right.$$
\indent
In the first subcase of the above we have the necessary estimate because
of the inequality (\arabic{FirstIndices}) and condition
$l({\cal B}) \geq 1$. In the second subcase we have the necessary estimate
from the inequality (\arabic{SecondIndices}) and again condition
$l({\cal B}) \geq 1$.\qed\\
\indent
This completes the proof of proposition~\ref{REstimate}.\\

\section{Compactness of ${\bold H}_r$.}\label{HCompactness}

\indent
From the definition of operator ${\bold H}_r$ we conclude that in order
to prove its compactness it suffices to prove compactness of each of the
terms below
$$\bar\partial_{\bold M} \vartheta^{\prime}_{\iota}(z)
\wedge R_r^{\iota}(\vartheta_{\iota}g)(z),\hspace{0.1in}
\vartheta^{\prime}_{\iota}(z) \cdot
R_{r+1}^{\iota}(\bar\partial_{\bold M} \vartheta_{\iota} \wedge g)(z)
\hspace{0.1in}\mbox{and}\hspace{0.1in}
\vartheta^{\prime}_{\iota}(z) \cdot H_r^{\iota}(\vartheta_{\iota}g)(z).$$
\indent
Compactness of the first two of these terms follows from the boundedness
of operators $R_r$ proved in proposition~\ref{REstimate} and compactness
of the embedding
$$\Gamma^{p, \alpha}(U) \rightarrow \Gamma^{p, \beta}(U)$$
for $\alpha > \beta$ \cite{Ad}. Another application of compactness of
this embedding shows that compactness of the third term
follows from the proposition below.

\begin{proposition}\label{HEstimate}
Let $0<\alpha< 1$, ${\bold M} \subset {\bold G}$ be a $C^{\infty}$ regular
q-concave CR submanifold of the form (\arabic{manifold}) and let $g \in \Gamma^{p, \alpha}_{(0,r)}({\bold M})$ be a form with compact support in
$U = {\widetilde U} \cap {\bold M}$.\\
\indent
Then for $r < q$ the operator $H_r$, defined in
(\arabic{LocalFormula}) satisfies the following estimate
$$\parallel H_r(g) \parallel_{\Gamma^{p,\alpha+1}_{(0,r)}(U)} <
C \cdot \parallel g \parallel_{\Gamma^{p, \alpha}_{(0,r)}(U)}
\eqno(\arabic{equation})
\newcounter{Hestimate}
\setcounter{Hestimate}{\value{equation}}
\addtocounter{equation}{1}$$
with a constant $C$ independent of $g$.
\end{proposition}

\indent
{\bf Proof.}\\
\indent
In our proof of proposition~\ref{HEstimate} we will
use the approximation of $H_r$ by the operators
$$H_r(\epsilon)(g)(z) = (-1)^{r} \frac{(n-1)!}{(2\pi i)^n}
\cdot \mbox{pr}_{\bold M} \circ
\int_{{\bold M}_{\epsilon}} \vartheta(\zeta) {\widetilde g}(\zeta)
\wedge \omega^{\prime}_{r} \left( \frac{P(\zeta,z)}
{\Phi(\zeta,z)}\right) \wedge\omega(\zeta)$$
\noindent
when $\epsilon$ goes to $0$.\\
\indent
Kernel of the operator $H_r(\epsilon)$ with the use of 
equalities (\arabic{dzetaF}) and (\arabic{dbarFbar}) may be represented
on ${\widetilde U} \times U$ as
$$\left. \vartheta(\zeta) \cdot
\omega^{\prime}_{r}\left(\frac{P(\zeta,z)}
{\Phi(\zeta,z)}\right) \wedge\omega(\zeta)
\right|_{{\widetilde U} \times U }
\eqno(\arabic{equation})
\newcounter{cauchykernel}
\setcounter{cauchykernel}{\value{equation}}
\addtocounter{equation}{1}$$
$$= \sum_{i,J} a_{(i,J)}(\zeta,z) \wedge
{\widetilde \phi}^{i,J}_{r}(\zeta, z) +
\sum_{i,J} b_{(i,J)}(\zeta,z) \wedge
{\widetilde \psi}^{i,J}_{r}(\zeta, z),$$
\noindent
where $i$ is an index, $J= \cup_{i=1}^{8} J_i$ is a
multiindex such that $i \not \in J,$
$a_{(i,J)}(\zeta,z)$ and $b_{(i,J)}(\zeta,z)$ are smooth functions of
$z$, $\zeta$ and $\theta(\zeta)$, and ${\widetilde \phi}^{i,J}_{r}(\zeta, z)$
and ${\widetilde \psi}^{i,J}_{r}(\zeta, z)$ are defined as follows:
$${\widetilde \phi}^{i,J}_{r}(\zeta, z) =\frac{1}
{{\Phi(\zeta,z)}^{n}} \times \sum \mbox{Det} \left[
Q^{(i)} {\overline F}^{(i)},\hspace{0.03in}
\overbrace{Q^{(j)} d_{\zeta}{\overline F}^{(j)}}^{j \in J_1},
\hspace{0.03in}
\overbrace{{\cal A}{\bar A} \cdot \bar\partial_{\zeta}a}^{j \in J_2},\right.
\eqno(\arabic{equation})
\newcounter{phiDet}
\setcounter{phiDet}{\value{equation}}
\addtocounter{equation}{1}$$
$$\left.\overbrace{a \cdot \mu_{\nu}}^{j \in J_3},\hspace{0.03in}
\overbrace{a \cdot \mu_{\tau}}^{j \in J_4}, \hspace{0.03in}
\overbrace{{\cal A}{\bar A} \cdot \bar \partial_z a}^{j \in J_5},
\hspace{0.03in}
\overbrace{a \cdot \bar\partial_z ({\cal A}{\bar A})}^{j \in J_6},
\hspace{0.03in}
\overbrace{{\overline F} \cdot \bar\partial_z Q}^{j \in J_7},
\hspace{0.03in}
\overbrace{Q \cdot \kappa}^{j \in J_8}
\right] \wedge\omega(\zeta),$$
\noindent
and
$${\widetilde \psi}^{i,J}_{r}(\zeta, z) =\frac{1}
{{\Phi(\zeta,z)}^{n}} \times \sum \mbox{Det} \left[
a_i {\cal A}{\bar A_i},\hspace{0.03in}
\overbrace{Q^{(j)} d_{\zeta}{\overline F}^{(j)}}^{j \in J_1},
\hspace{0.03in}
\overbrace{{\cal A}{\bar A} \cdot \bar\partial_{\zeta}a}^{j \in J_2},\right.
\eqno(\arabic{equation})
\newcounter{psiDet}
\setcounter{psiDet}{\value{equation}}
\addtocounter{equation}{1}$$
$$\left.\overbrace{a \cdot \mu_{\nu}}^{j \in J_3},\hspace{0.03in}
\overbrace{a \cdot \mu_{\tau}}^{j \in J_4}, \hspace{0.03in}
\overbrace{{\cal A}{\bar A} \cdot \bar \partial_z a}^{j \in J_5},
\hspace{0.03in}
\overbrace{a \cdot \bar\partial_z ({\cal A}{\bar A})}^{j \in J_6},
\hspace{0.03in}
\overbrace{{\overline F} \cdot \bar\partial_z Q}^{j \in J_7},
\hspace{0.03in}
\overbrace{Q \cdot \kappa}^{j \in J_8}
\right] \wedge\omega(\zeta).$$
\indent
Applying then (\arabic{dzetaF}) we conclude that
${\widetilde \phi}^{i,J}_r(\zeta, z)$
can be represented as a finite sum of the following terms with smooth
coefficients in $\zeta, z$ and $\theta(\zeta)$:
$$\phi^{i,J}_r(\zeta, z) =\frac{{\overline F}^{(i)}}
{{\Phi(\zeta,z)}^n} \cdot
\overbrace{d_{\zeta}{\overline F}^{(j)}}^{j \in J_1} \wedge
\overbrace{d_{\zeta}F^{(j)}}^{j \in J_1} \wedge
\wedge \overbrace{{\cal A}{\bar A} \cdot
\bar\partial_{\zeta} a}^{j \in J_2} \wedge
\overbrace{a \cdot \mu_{\nu}}^{j \in J_3}\wedge
\overbrace{a \cdot \mu_{\tau}}^{j \in J_4}
\eqno(\arabic{equation})
\newcounter{phiForm}
\setcounter{phiForm}{\value{equation}}
\addtocounter{equation}{1}$$
$$\wedge \overbrace{{\cal A}{\bar A} \cdot \bar \partial_z a}^{j \in J_5}
\wedge\overbrace{a \cdot \bar\partial_z ({\cal A}{\bar A})}^{j \in J_6}
\wedge \overbrace{{\overline F} \cdot \bar\partial_z Q}^{j \in J_7}
\wedge \overbrace{\kappa}^{j \in J_8}\wedge
\overbrace{d\zeta}^{n-|J_1|}.$$\\
\indent
Similarly, ${\widetilde \psi}^{i,J}_r(\zeta, z)$ can be represented as a
finite sum of the following terms with smooth coefficients in $\zeta, z$
and $\theta(\zeta)$:
$$\psi^{i,J}_r(\zeta, z) =\frac{a_i {\cal A}{\bar A_i}}
{{\Phi(\zeta,z)}^n} \cdot
\overbrace{d_{\zeta}{\overline F}^{(j)}}^{j \in J_1} \wedge
\overbrace{d_{\zeta}F^{(j)}}^{j \in J_1} \wedge
\wedge \overbrace{{\cal A}{\bar A} \cdot
\bar\partial_{\zeta} a}^{j \in J_2} \wedge
\overbrace{a \cdot \mu_{\nu}}^{j \in J_3}\wedge
\overbrace{a \cdot \mu_{\tau}}^{j \in J_4}
\eqno(\arabic{equation})
\newcounter{psiForm}
\setcounter{psiForm}{\value{equation}}
\addtocounter{equation}{1}$$
$$\wedge \overbrace{{\cal A}{\bar A} \cdot \bar \partial_z a}^{j \in J_5}
\wedge\overbrace{a \cdot \bar\partial_z ({\cal A}{\bar A})}^{j \in J_6}
\wedge \overbrace{{\overline F} \cdot \bar\partial_z Q}^{j \in J_7}
\wedge \overbrace{\kappa}^{j \in J_8}\wedge
\overbrace{d\zeta}^{n-|J_1|}.$$\\
\indent
As in the proof of proposition~\ref{REstimate} we use estimates
(\arabic{termsestimates}) for the terms of the determinants in
(\arabic{phiDet}) and (\arabic{psiDet}) and obtain the representations
$$a_{(i,J)}(\zeta,z) \wedge {\widetilde g}(\zeta) \wedge
\phi^{i,J}_{r}(\zeta, z)$$
$$= \sum_{|T|+|E| \leq |J_8|+1}
c_{ \{I,d,j \} }(\zeta, z) {\widetilde g}(\zeta)
\{\mbox{Im}F(\zeta, z) \}^{T} {\cal K}^{I(J,S,T,E)}_{d,j}(\zeta, z),
\eqno(\arabic{equation})
\newcounter{phiRepresentation}
\setcounter{phiRepresentation}{\value{equation}}
\addtocounter{equation}{1}$$
and
$$b_{(i,J)}(\zeta,z) \wedge {\widetilde g}(\zeta) \wedge
\psi^{i,J}_{r}(\zeta, z)$$
$$= \sum_{|T|+|E| \leq |J_8|}
c_{ \{I,d,j \} }(\zeta, z) {\widetilde g}(\zeta)
\{\mbox{Im}F(\zeta, z) \}^{T} {\cal K}^{I(J,S,T,E)}_{d,j}(\zeta, z).
\eqno(\arabic{equation})
\newcounter{psiRepresentation}
\setcounter{psiRepresentation}{\value{equation}}
\addtocounter{equation}{1}$$
\indent
Multiindices $T$ and $E$ in (\arabic{phiRepresentation}) are
obtained from the decomposition
$$\{ \overline{F}(\zeta, z) \}^{ \{J_7 \cup i \} }$$
$$= \sum_{|T|+|E|+\frac{1}{2}\left(|G|+|H|\right)=|J_7|+1}
c_{\left\{ T,E,G,H \right\} }(\zeta, z)
\{ \mbox{Im}F(\zeta,z)\}^{T} \{ \rho(\zeta) \}^{E}
(\zeta - z)^{G} (\bar{\zeta} - \bar{z})^{H}$$
and multiindices $I_i$ for $i= 1, \dots, 5$ and indices $d,j$ in
(\arabic{phiRepresentation}) satisfy the conditions below
$$\begin{array}{lllllll}
\hspace{0.3in}d = 0,\vspace{0.1in}\\
\hspace{0.3in}j = n,\vspace{0.1in}\\
\hspace{0.3in}|I_1| = |E|,\vspace{0.1in}\\
\hspace{0.3in}|I_2| + |I_3| = 3|J_2|+3|J_3|+2|J_4|+3|J_5|+2|J_6|+|J_8|
\vspace{0.1in}\\
+2\left( |J_7|+1-|T|-|E| \right),\vspace{0.1in}\\
\hspace{0.3in}|I_4| = |J_1|,\vspace{0.1in}\\
\hspace{0.3in}|I_5| = 0.
\end{array}
\eqno(\arabic{equation})
\newcounter{KphiIndices}
\setcounter{KphiIndices}{\value{equation}}
\addtocounter{equation}{1}$$
\indent
Multiindices $T$ and $E$ in (\arabic{psiRepresentation}) are
obtained from the decomposition
$$\{ \overline{F}(\zeta, z) \}^{ \{J_7 \} }$$
$$= \sum_{|T|+|E|+\frac{1}{2}\left(|G|+|H|\right)=|J_7|}
c_{\left\{ T,E,G,H \right\} }(\zeta, z)
\{ \mbox{Im}F(\zeta,z)\}^{T} \{ \rho(\zeta) \}^{E}
(\zeta - z)^{G} (\bar{\zeta} - \bar{z})^{H}$$
and multiindices $I_i$ for $i= 1, \dots, 5$ and indices $d,j$ in
(\arabic{psiRepresentation}) satisfy the conditions
$$\begin{array}{lllllll}
\hspace{0.3in}d = 0,\vspace{0.1in}\\
\hspace{0.3in}j = n,\vspace{0.1in}\\
\hspace{0.3in}|I_1| = |E|,\vspace{0.1in}\\
\hspace{0.3in}|I_2| + |I_3| =
3+3|J_2|+3|J_3|+2|J_4|+3|J_5|+2|J_6|+|J_8|
\vspace{0.1in}\\
+2\left( |J_7|-|T|-|E| \right),\vspace{0.1in}\\
\hspace{0.3in}|I_4| = |J_1|,\vspace{0.1in}\\
\hspace{0.3in}|I_5| = 0.
\end{array}
\eqno(\arabic{equation})
\newcounter{KpsiIndices}
\setcounter{KpsiIndices}{\value{equation}}
\addtocounter{equation}{1}$$
\indent
Using representations (\arabic{phiRepresentation}) and
(\arabic{psiRepresentation}) we reduce the statement
of proposition~\ref{HEstimate} to the same statement for each term
$${\widetilde g}(\zeta){\cal B}^{T,I}_{d,j}(\zeta, z)$$
of the right hand side of these representations.\\
\indent
Proceeding then as in the proof of proposition~\ref{REstimate} we reduce the
problem to lemma~\ref{GammaAlpha}. To check applicability of
lemma~\ref{GammaAlpha} we have to prove only the second condition of
(\arabic{FinalIndices}) for the kernels obtained from
$\phi^{i,J}_{r}(\zeta, z)$ and $\psi^{i,J}_{r}(\zeta, z)$ and satisfying
(\arabic{KphiIndices}) and (\arabic{KpsiIndices}) respectively. The
reason for that is that for these kernels we have
$k({\cal B}) \leq 0 \leq 2n-2m-1$ and therefore only the second condition
of (\arabic{FinalIndices}) is needed for application of
lemma~\ref{Integral}.\\
\indent
For multiindices satisfying (\arabic{KphiIndices}) we have
$$2n-k({\cal B})-2h({\cal B})+2l({\cal B})+2s({\cal B})$$
$$=2n+|I_2|+|I_3|-4n+2|T|+2|I_1|+2|I_4|$$
$$=|J_2|+|J_3|+|J_5|-|J_8| \geq |J_1|+|J_2|+|J_3|-m+1,$$
\noindent
where we used the inequality
$$|J_1|+|J_8| \leq m-1.$$
\indent
Therefore, if
$$2n-k({\cal B})-2h({\cal B})+2l({\cal B})+2s({\cal B}) \leq 0$$
\noindent
then
$$|J_1|+|J_2|+|J_3| \leq m-1$$
\noindent
and hence
$$|J_5| \geq n-r-1-|J_1|-|J_2|-|J_3| \geq n-r-m > n-q-m,$$
\noindent
which is impossible.\\
\indent
Analogously, for multiindices satisfying (\arabic{KpsiIndices}) we have
$$2n-k({\cal B})-2h({\cal B})+2l({\cal B})+2s({\cal B})$$
$$=2n+|I_2|+|I_3|-4n+2|T|+2|I_1|+2|I_4|$$
$$=|J_2|+|J_3|+|J_5|-|J_8|+1 \geq |J_1|+|J_2|+|J_3|-m+1,$$
\noindent
where we used the inequality
$$|J_1|+|J_8| \leq m.$$
\indent
The same arguments as above again show that
$$2n-k({\cal B})-2h({\cal B})+2l({\cal B})+2s({\cal B}) \geq 1.$$
\qed


\begin{thebibliography}{width}
\bibitem[Ad]{Ad} R. A. Adams, Sobolev Spaces, Academic Press, New York, 1975.
\bibitem[AH]{AH} A. Andreotti, C.D. Hill, E. E. Levi convexity
and the Hans Lewy problem, I, II, Ann. Scuola Norm. Super. Pisa,
26:2 (1972), 325-363, 26:4 (1972), 747-806.
\bibitem[AiH]{AiH} R.A. Airapetian, G.M. Henkin, Integral 
representation of differential forms on Cauchy-Riemann manifolds
and the theory of CR functions, I, II, Russ. Math. Surv., 39 (1984),
41-118, Math.USSR Sbornik 55:1 (1986), 91-111.
\bibitem[AnH]{AnH} A. Andreotti, C.D. Hill, E. E. Levi convexity
and the Hans Lewy problem, I, II, Ann. Scuola Norm. Super. Pisa,
26:2 (1972), 325-363, 26:4 (1972), 747-806.
\bibitem[Ba]{Ba} M. Y. Barkatou, Optimal Regularity for $\bar\partial_b$ on
CR manifolds, Preprint, 1995.
\bibitem[BGG]{BGG} R. Beals, B. Gaveau, P.C. Greiner, 
The Green function of model step two hypoelliptic operators and the analysis
of certain Cauchy-Riemann complexes, Preprint, 1995.
\bibitem[B]{B} A. Boggess, CR Manifolds and the Tangential
Cauchy-Riemann Complex, CRC Press, Boca Raton, Florida, 1991.
\bibitem[FS]{FS} G.B. Folland, E.M. Stein, Estimates for
the $\bar\partial_b$ complex and analysis on the Heisenberg group,
Comm. Pure Appl. Math. 27 (1974), 429-522.
\bibitem[H1]{H1} G.M. Henkin, The Hans Lewy equation and
analysis on pseudoconvex manifolds, Math.USSR Sbornik 31 (1977), 59-130.
\bibitem[H2]{H2} G.M. Henkin, Solution des equations de
Cauchy-Riemann tangentielle sur les varietes de Cauchy-Riemann
q-concaves, C. R. Acad. Sci. Paris, 292 (1981), 27-30.
\bibitem[H3]{H3} G.M. Henkin, The Method of Integral
Representations in Complex Analysis, in Encyclopedia of Mathematical
Sciences, volume 7, Several Complex Variables I, Springer-Verlag, 1990.
\bibitem[J]{J} H. Jacobowitz, An Introduction to CR Structures,
Math. Surveys and Monographs 32, AMS, Providence, Rhode Island.
\bibitem[K]{K} J.J. Kohn, Boundaries of complex manifolds, Proc.
Conf. Complex Analysis, Minneapolis 1964, Berlin-Heidelberg-New-York,
Springer Verlag, 1965, 81-94.
\bibitem[KR]{KR} J.J. Kohn, H. Rossi, On the extension of
holomorphic functions from the boundary of a complex manifold, Ann. of
Math. 81 (1965), 451-472.
\bibitem[Na]{Na} I. Naruki, Localization principles for differential
complexes and its applications, Publ. Res. Inst. Math., Kyoto, ser. A,
8:1 (1972), 43-110.
\bibitem[P]{P} P. Polyakov, Sharp estimates for
operator $\bar\partial_{\bold M}$ on a q-concave CR manifold, to appear
in The Journal of Geometric Analysis.
\bibitem[Siu]{Siu} Y.-T. Siu, The $\bar\partial$ problem with uniform bounds
on derivatives, Math.Ann. 207 (1974), 163-176.
\bibitem[S]{S} E.M. Stein, Singular integrals and estimates for
the Cauchy-Riemann equations, Bull.Amer.Math.Soc. 79:2 (1973),
440-445.
\end{thebibliography}
\end{document}